\documentclass[12pt,psamsfonts,leqno,oneside,letterpaper]{amsart}
\usepackage[dvips,text={6.5truein,9truein},left=1truein,top=1truein]{geometry}
\usepackage{amssymb,amsmath,amscd,enumerate,url}

\usepackage[pdftex]{graphicx}
\usepackage[colorlinks,linkcolor=blue,citecolor=blue,pdfstartview=FitH]{hyperref}
\input xy
\xyoption{all}
\SelectTips{cm}{12}

\title{Volume and homology of one-cusped hyperbolic $3$-manifolds}

\author{Marc Culler} 
\address{Department of Mathematics, Statistics, and Computer Science (M/C 249)\\
University of Illinois at Chicago\\
851 S. Morgan St.\\
Chicago, IL 60607-7045}
\email{culler@math.uic.edu}
\urladdr{http://www.math.uic.edu/~culler}
\author{Peter B. Shalen}
\address{Department of Mathematics, Statistics, and Computer Science (M/C 249)\\
University of Illinois at Chicago\\
851 S. Morgan St.\\
Chicago, IL 60607-7045}
\email{shalen@math.uic.edu}
\urladdr{http://www.math.uic.edu/~shalen}

\theoremstyle{definition}
\newtheorem{para}{}[section]

\newtheorem{remark}[para]{Remark}
\newtheorem{remarks}[para]{Remarks}
\newtheorem{notation}[para]{Notation}
\newtheorem{convention}[para]{Convention}
\newtheorem{definition}[para]{Definition}
\newtheorem{definitions}[para]{Definitions}

\theoremstyle{plain}
\newtheorem{theorem}[para]{Theorem}
\newtheorem{lemma}[para]{Lemma}
\newtheorem{proposition}[para]{Proposition}
\newtheorem{corollary}[para]{Corollary}
\newtheorem{conjecture}[para]{Conjecture}
\numberwithin{equation}{para}
\newtheorem{claim}[equation]{}

\numberwithin{equation}{para}

\newenvironment{conclusions}{

\begin{enumerate}
}{\end{enumerate}}

\newenvironment{conditions}{

\begin{enumerate}
}{\end{enumerate}}

\newenvironment{alternatives}{

\begin{enumerate}
}{\end{enumerate}}

\newcommand\Number{\begin{para}}
\newcommand\EndNumber{\end{para}}
\newcommand\Definition{\begin{definition}}
\newcommand\EndDefinition{\end{definition}}
\newcommand\Definitions{\begin{definitions}}
\newcommand\EndDefinitions{\end{definitions}}
\newcommand\Theorem{\begin{theorem}}
\newcommand\EndTheorem{\end{theorem}}
\newcommand\Conjecture{\begin{conjecture}}
\newcommand\EndConjecture{\end{conjecture}}
\newcommand\Remark{\begin{remark}}
\newcommand\EndRemark{\end{remark}}
\newcommand\Remarks{\begin{remarks}}
\newcommand\EndRemarks{\end{remarks}}
\newcommand\Convention{\begin{convention}}
\newcommand\EndConvention{\end{convention}}
\newcommand\Notation{\begin{notation}}
\newcommand\EndNotation{\end{notation}}
\newcommand\Lemma{\begin{lemma}}
\newcommand\EndLemma{\end{lemma}}
\newcommand\Proposition{\begin{proposition}}
\newcommand\EndProposition{\end{proposition}}
\newcommand\Corollary{\begin{corollary}}
\newcommand\EndCorollary{\end{corollary}}
\newcommand\Claim{\begin{claim}}
\newcommand\EndClaim{\end{claim}}
\newcommand\Proof{\begin{proof}}
\newcommand\EndProof{\end{proof}}
\newcommand\Equation{\begin{equation}}
\newcommand\EndEquation{\end{equation}}

\newcommand\Conclusions{\begin{conclusions}}
\newcommand\EndConclusions{\end{conclusions}}
\newcommand\Alternatives{\begin{alternatives}}

\newcommand\twoareminusthree{6.3}
\newcommand\exactsequence{6.1}
\newcommand\geesubem{6.2}
\newcommand\topeleven{8.13}

\newcommand\EndAlternatives{\end{alternatives}}
\newcommand\Conditions{\begin{conditions}}
\newcommand\EndConditions{\end{conditions}}
\newcommand\newrk{\dim_{\ZZ_2}}
\newcommand\mfldrk{{\rm rk}_2}
\newcommand\Letters{\begin{letters}}
\newcommand\EndLetters{\end{letters}}
\newcommand\Bullets{\begin{itemize}}
\newcommand\EndBullets{\end{itemize}}

\renewcommand\Im{\mathop{\rm Im}}

\newcommand\whatvolume{5.06}
\newcommand\inter{\mathop{\rm int}}

\newcommand\Mthick{M_{\rm thick}}
\newcommand\Mthin{M_{\rm thin}}

\newcommand\tM{\widetilde M}
\newcommand\tP{\widetilde P}
\newcommand\tGamma{\widetilde \Gamma}
\newcommand\tT{\widetilde T}
\newcommand\tN{\widetilde N}

\newcommand\tMthick{\widetilde M_{\rm thick}}
\newcommand\tMthin{\widetilde M_{\rm thin}}

\newcommand\calc{{\mathcal C}}

\newcommand\Isom{{\rm Isom}}
\newcommand\arccosh{{\rm arccosh}}

\newcommand\arctanh{{\rm arctanh}}
\newcommand\arcsech{{\rm arcsech}}
\newcommand\arcsinh{{\rm arcsinh}}
\newcommand\sech{{\rm sech}}
\newcommand\arcsec{{\rm arcsec}}
\newcommand\ZZ{{\mathbb Z}}
\newcommand\NN{{\mathbb N}}
\newcommand\RR{{\mathbb R}}
\newcommand\CC{{\mathbb C}}

\newcommand\HH{{\mathbb H}}

\newcommand\dist{\mathop{\rm dist}}
\newcommand\genus{\mathop{\rm genus}}
\newcommand\vol{\mathop{\rm vol}}

\newcommand\inj{{\rm inj}}

\newcommand\length{\mathop{{\rm length}}}
\newcommand\rank{\mathop{{\rm rank}}}

\newcommand\trace{\mathop{\rm trace}}

\begin{document}

\begin{abstract}
Let $M$ be a complete, finite-volume, orientable hyperbolic manifold
having exactly one cusp. If we assume that $\pi_1(M)$ has no subgroup
isomorphic to a genus-$2$ surface group, and that either (a)
$\dim_{\ZZ_p}H_1(M;\ZZ_p)\ge 5$ for some prime $p$, or (b)
$\dim_{\ZZ_2}H_1(M;\ZZ_2)\ge 4$, and the subspace of $H^2(M;\ZZ_2)$
spanned by the image of the cup product $H^1(M;\ZZ_2)\times
H^1(M;\ZZ_2)\to H^2(M;\ZZ_2)$ has dimension at most $1$, then $\vol
M>\whatvolume.$ If we assume that $\dim_{\ZZ_2}H_1(M;\ZZ_2)\ge 7$, and
that the compact core $N$ of $M$ does not contain a genus-$2$ closed
incompressible surface, then $\vol M>\whatvolume.$
\end{abstract}

\maketitle

\section{Introduction}

In this paper we shall prove:

\Theorem\label{motivating result}
Let $M$ be a complete, finite-volume, orientable hyperbolic
$3$-manifold having exactly one cusp, such that that $\pi_1(N)$ has no
subgroup isomorphic to a genus-$2$ surface group. Suppose that either
\Alternatives
\item $\dim_{\ZZ_p}H_1(N;\ZZ_p)\ge 5$ for some prime $p$, or
\item $\dim_{\ZZ_2}H_1(N;\ZZ_2)\ge 4$, and the subspace of
  $H^2(M;\ZZ_2)$ spanned by the image of the cup product
  $H^1(M;\ZZ_2)\times H^1(M;\ZZ_2)\to H^2(M;\ZZ_2)$ has dimension at
  most $1$.
\EndAlternatives
Then
$$\vol M>\whatvolume.$$ \EndTheorem

Contrapositively, one can think of Theorem \ref{motivating
result} as saying that if the one-cusped hyperbolic $3$-manifold has
volume at most $\whatvolume$, then $\dim_{\ZZ_p}H_1(N;\ZZ_p)\le 4$ for every
prime $p$; and that if in addition one assumes that the span of the
image of the cup product has dimension at most $1$, then
$\dim_{\ZZ_2}H_1(N;\ZZ_2)\le 3$.

We do not know whether these bounds on homology are sharp. The
Weeks-Hodgson census \cite{snappea} contains many examples of
one-cusped manifolds with volume $<\whatvolume$ for which
$\dim_{\ZZ_2}H_1(N;\ZZ_2)= 3$. We have not calculated the cup product
for these examples.

Theorem \ref{motivating result} is an immediate consequence of the
following two results, Proposition \ref{k-free filling} and Theorem
\ref{main estimate}.

Recall that the {\it rank} of a group $\Phi$ is the smallest
cardinality of any generating set of $\Phi$. A group $\Gamma$ is
said to be {\it $k$-free}\/, where $k$ is a given integer, if every
subgroup of $\Gamma$ whose rank is at most $k$ is a free group.

\Proposition\label{k-free filling}
Let $N$ be a compact, orientable $3$-manifold whose boundary is a
torus, and let $k\ge2$ be an integer. Suppose that $\pi_1(N)$ has no
subgroup isomorphic to a genus-$g$ surface group for any integer
$g$ with $k/2<g<k$. In addition, suppose  that     either
\Alternatives
\item $\dim_{\ZZ_p}H_1(N;\ZZ_p)\ge k+2$ for some prime $p$, or
\item $\dim_{\ZZ_2}H_1(N;\ZZ_2)\ge k+1$, and the subspace of
  $H^2(N;\ZZ_2)$ spanned by the image of the cup product
  $H^1(M;\ZZ_2)\times H^1(M;\ZZ_2)\to H^2(M;\ZZ_2)$ has dimension at
  most $k-2$.
\EndAlternatives

Then there is an infinite sequence of closed manifolds
$(M_j)_{j\in\NN}$, obtained by distinct Dehn fillings of $N$, such
that $\pi_1(M_j)$ is $k$-free for every $j\ge1$.
\EndProposition

\Theorem\label{main estimate}
Let $M$ be a complete, finite-volume, orientable hyperbolic manifold
having exactly one cusp. Suppose that there is an infinite
sequence of closed manifolds $(M_j)_{j\in\NN}$, obtained by distinct
Dehn fillings of  the compact core of $M$, such that $\pi_1(M_j)$ is
$3$-free for every $j\ge0$. Then 
$$\vol M>\whatvolume.$$ \EndTheorem

The proof of Proposition \ref{k-free filling} will be given in Section
\ref{k-free filling section}. The proof of Theorem
\ref{main estimate} will be given in Section \ref{main estimate section}.

In Section \ref{singular applications} we prove a result, Proposition
\ref{lots of H_1}, which asserts that if a one-cusped (complete,
finite-volume, orientable) hyperbolic manifold $M$ satisfies
$\dim_{\ZZ_2}H_1(M;\ZZ_2)\ge 7$, then either $\vol M>\whatvolume$ or
$N$ contains a genus-$2$ (embedded) incompressible surface.  We may
think of this as saying that the restriction on surface subgroups in
Theorem \ref{motivating result} can be relaxed at the expense of
strengthening the lower bound on homology.  The proof combines Theorem
\ref{main estimate} with the results of \cite{last}.

In \cite{cds}, the following strictly stronger version of Proposition
\ref{lots of H_1} will be proved: If $M$ is a complete, finite-volume,
orientable hyperbolic manifold having exactly one cusp, such that
$\dim_{\ZZ_2}H_1(M;\ZZ_2)\ge 7$, then $\vol M>\whatvolume$.  The proof
of the stronger result takes Proposition \ref{lots of H_1} as a
starting point, and involves a careful analysis of the case where $M$
contains an incompressible surface of genus $2$.

Our proof of Proposition \ref{k-free filling} depends on the following
result, Theorem \ref{dewey}, which is of independent interest and
is related to a result of Oertel \cite[Proposition 1.1]{oertel}.

\Definition\label{singular boundary slope def} Let $N$ be a compact,
orientable $3$-manifold whose boundary is a torus, and let $g>0$ be an
integer. A slope (see \ref{slope def}) $\alpha$ in $\partial N$ will be
called a {\it genus-$g$ singular boundary slope} if there exist a
compact orientable surface $F$, with genus $g$ and non-empty boundary,
and a map $f:F\to N$, such that (i) $f_\sharp:\pi_1(F)\to\pi_1(N)$ and
$f_\sharp:\pi_1(F,\partial F)\to\pi_1(N,\partial N)$ are injective,
and (ii) $f$ maps each component of $\partial F$ homeomorphically onto
a simple closed curve in $\partial N$ representing the slope
$\alpha$.  \EndDefinition

(We show in Proposition \ref{same difference} that injectivity of
$f_\sharp:\pi_1(F)\to\pi_1(N)$ implies  injectivity of $f_\sharp:\pi_1(F\partial
F)\to\pi_1(N,\partial N)$.)

\Theorem\label{dewey} Let $N$ be a compact, orientable $3$-manifold
whose boundary is a torus, let $\alpha$ be a slope (see \ref{slope
  def}) in $\partial N$, and let $g>0$ be an integer.  Suppose that
the fundamental group of the manifold $N(\alpha)$ obtained from $N$ by
the surgery corresponding to $\alpha$ (see \ref{slope def}) has a
subgroup isomorphic to a genus-$g$ surface group. Then either
$\pi_1(N)$ has a subgroup isomorphic to a genus-$g$ surface group, or
$\alpha$ is a singular genus-$g$ boundary slope.  \EndTheorem

The proof of Theorem \ref{dewey} will also be given in Section
\ref{k-free filling section}.

A crucial ingredient in the proof of Theorem \ref{main estimate}
is a result, Proposition \ref{popeye}, which asserts that if $c$ is a
suitably short geodesic in a hyperbolic $3$-manifold $M$, and
$\pi_1(M)$ is $k$-free for a given $k\ge3$, then there exists a point
$Q\in M$ which is  far from $c$; the precise lower bound for
the distance from $Q$ to (a suitably chosen point of) $c$ is a
function of $k$ and of the length $l$ of $c$, which increases when one
decreases  $l$ or increases $k$. For small $l$, this function is
strictly greater than the lower bound for the radius of an embedded
tube about $c$ given by \cite[Corollary 10.2]{accs}. Proposition
\ref{popeye} also asserts that $Q$ may be taken to be a ($\log3)$-thick point,
in the sense that it is the center of an embedded hyperbolic ball of
radius $(\log3)/2$. 

The proof of Proposition \ref{popeye} is based on a new application of
the Generalized $\log(2k-1)$ Theorem (\cite[Theorem 4.1]{accs}, which
we combine with the Tameness Theorem of \cite{agol} or \cite{cg}), and
may be regarded as an analogue of \cite[Corollary 10.2]{accs} for $k>2$.

We prove Proposition \ref{popeye} in Section \ref{distant section},
after some algebraic preliminaries in Section \ref{algebra
section}. We deduce Proposition \ref{popeye}, via the topology of
Margulis tubes, from a result, Proposition \ref{big radius}, which
applies in a more general setting and which we expect to have other
applications.

In Section \ref{moe}, we combine Proposition \ref{popeye} with
geometric convergence techniques based on Thurston's Dehn filling
theorem to prove a result, Theorem \ref{geometric fact}, which asserts
that if $M$ satisfies the hypotheses of Theorem \ref{main estimate},
there is a ($\log3)$-thick point $Q$ which is suitably distant from a
maximal standard neighborhood  (\ref{distance notions}) of the cusp of
$M$. Much of the general information about Dehn filling and geometric
convergence that we needed is contained in our Proposition \ref{tube
convergence}.

In Section \ref{main estimate section} we deduce Theorem \ref{main
estimate} from Theorem \ref{geometric fact}. This step requires the
use of a stronger form of B\"or\"ockzy's results on density of
sphere-packing \cite{bor} than has previously been used in the study
of volumes of hyperbolic $3$-manifolds. The needed background is
explained in Section \ref{packing section}.

Section \ref{convention section} is devoted to establishing some
conventions that are used in the rest of the paper.  In section
\ref{cup product section} we establish a variant of \cite[Lemma
\twoareminusthree]{last} which is needed for the proof of case (b) of
Theorem \ref{motivating result}.

We are grateful to Rosemary Guzman for correcting a number of
errors in an earlier version of the paper.  This research was
partially supported by NSF grants DMS-0204142 and DMS-0504975.

\section{General conventions}\label{convention section}

\Number\label{topology notions}
A {\it manifold} may have a boundary. If $M$ is a manifold, we shall
denote the boundary of $M$ by $\partial M$ and its interior
$M-\partial M$ by $\inter M$.

In some of our topological results about manifolds of dimension
$\le3$, specifically those in Section \ref{k-free filling section}, we do not
specify a category. These results may be interpreted in the category
in which manifolds are topological, PL or smooth, and submanifolds are
respectively locally flat, PL or smooth; these three categories are
equivalent in these low dimensions as far as classification is
concerned. Some of the proofs in Section 3 are done in the PL
category, but the applications to hyperbolic manifolds in later sections
are carried out in the smooth category.
\EndNumber

\Definitions\label{simple def}
A $3$-manifold $M$ will be termed {\it irreducible} if every
$2$-sphere in $M$ bounds a ball in $M$.  We shall say that $M$ is {\it
boundary-irreducible} if $\partial M$ is $\pi_1$-injective in $M$, or
equivalently if for every properly embedded disk $D\subset M$ the
simple closed curve $\partial D$ bounds a disk in $\partial M$. We
shall say that a $3$-manifold $M$ is {\it simple} if (i) $M$ is
compact, connected, orientable, irreducible and boundary-irreducible;
(ii) every non-cyclic abelian subgroup of $\pi_1(M)$ is carried
(up to conjugacy) by a torus component of $\partial M$; and
(iii) $M$ is not a closed manifold with finite fundamental
group, and is not homeomorphic to $S^1\times S^1\times[0,1]$.
\EndDefinitions

\Number\label{slope def}
If $X$ is a $2$-torus, we define a {\it slope} in $X$ to be an
unoriented isotopy class of non-trivial simple closed curves in
$X$. 

If $N$ is a $3$-manifold, $X$ is a torus component of $\partial N$,
and $\alpha$ is a slope in $X$, we shall denote by $N(\alpha)$ the
manifold obtained from the disjoint union of $N$ and $D^2\times S^1$
by identifying $D^2\times S^1$ with $X$ via a homeomorphism that maps
$D^2\times\{1\}$ onto a curve with slope $\alpha$. We shall say that
$N(\alpha)$ is obtained from $N$ by the {\it Dehn filling
  corresponding to $\alpha$}.
\EndNumber

\Number\label{distance notions}If $P$ and $Q$ are points in a complete
hyperbolic manifold $M$, we shall denote by $\dist_M(P,Q)$ the
hyperbolic distance between $P$ and $Q$ (which may be thought of as
the length of a shortest geodesic from $P$ to $Q$, if $P\ne Q$).

Let $X$ be a subset of a complete hyperbolic $3$-manifold $M$. For any
point $P\in M$ we define $\dist(P,X)=\inf_{Q\in X}\dist(P,Q)$. We have
$\dist(P,X)=\dist(P,\bar X)$, and $\dist(P,X)=0$ if and only if
$P\in\bar X$. If $r$ is a real number, we shall denote by $N(r,X)$ the
set of all points $P\in M$ such that $\dist(P,X)<r$. Note that
$N(r,X)=\emptyset$ if $r\le0$.  If $Q$ is a point of $M$ and $r$ is a
real number, we shall write $N(r,Q)=N(r,\{Q\})$.

By a {\it cylinder} in $\HH^n$ we mean a set of the form $N(r,A)$,
where $r$ is a (strictly) positive real number and
$A$ is a line in $\HH^n$.

Let $M$ be a complete hyperbolic $n$-manifold. We shall often identify
$M$ with a quotient $\HH^n/\Gamma$ where $\Gamma\le\Isom_+(\HH^n)$ is
discrete and torsion-free.  We define a {\it tube} (or, respectively,
a {\it standard cusp neighborhood}) in $M$ to be the image in $M$ of a
cylinder (or, respectively, a horoball) in $\HH^n$ which is precisely
invariant under $\Gamma$. Note that tubes and standard cusp
neighborhoods are open subsets of $M$.

Every simple closed geodesic in $M$ has a neighborhood which is a
tube.

An orbit of points in the sphere at infinity which are fixed points of
parabolic elements of $\Gamma$ is said to define a {\it cusp}
$\mathcal K$ of $M$. If $z$ is a point of the given orbit, some
horoball based at $z$ is precisely invariant under $\Gamma$ and thus
defines a standard cusp neighborhood, which we may also refer to as a
standard neighborhood of $\mathcal K$.

The tubes about a given simple geodesic $c$, or the standard
neighborhoods of a given cusp $\mathcal K$, are totally ordered by
inclusion, and their union is the unique tube about $c$, or
respectively the unique maximal standard neighborhood of $\mathcal K$.

\EndNumber

\Number\label{cylinder notions} Let $M=\HH^n/\Gamma$ be a complete
hyperbolic $n$-manifold, and let $q:\HH^n\to M$ denote the quotient
map.

If $\lambda$ is a positive number, we shall denote by
$\Mthin(\lambda)$ the set of all points of $M$ at which the
injectivity radius is strictly less than $\lambda/2$.  We shall
denote by $\tMthin(\lambda)$ the set $q^{-1}(\Mthin(\lambda))$.
Equivalently, $\tMthin(\lambda)$ is the set of all points $z\in\HH^3$
such that $\dist(\gamma( z), z)<\lambda$ for some
$\gamma\in\Gamma-\{1\}$.
  
We set $\Mthick(\lambda)=M-\Mthin(\lambda)$ and
$\tMthick(\lambda)=q^{-1}(\Mthick(\lambda))=\HH^n -\tMthin(\lambda)$.

Each non-trivial element $\gamma$ of $\Gamma$ lies in a unique maximal
abelian subgroup, which is the centralizer of $\gamma$. In particular,
non-trivial elements which lie in distinct maximal abelian subgroups
do not commute.

The non-trivial elements of a maximal abelian subgroup $C$ either are
all loxodromic isometries with a common axis, which we shall denote by
$A_C$, or are all parabolic isometries with a common fixed point at
infinity, which we shall denote by $\omega_C$. If $C$ is loxodromic it
must be cyclic.  We shall denote by $\calc(\Gamma)=\calc(M)$ the set
of all maximal abelian subgroups of $\Gamma$.  For any positive real
number $\lambda$, we denote by
$\calc_\lambda(\Gamma)=\calc_\lambda(M)$ the set of all subgroups
$C\in\calc(\Gamma)$ which are either parabolic or are generated by
loxodromic elements of translation length $<\lambda$.

Let $\gamma$ be an element of $\Gamma$, and let $\lambda$ be a
positive real number.  We shall denote by $Z_\lambda(\gamma)$ the set
of points $ z\in\HH^3$ such that $\dist(z,\gamma(z))<\lambda$.  If
$\gamma$ is parabolic then $Z_\lambda(\gamma)$ is a horoball based at
the fixed point of $\gamma$. If $\gamma$ is loxodromic and the
translation length $l$ of $\gamma$ is at least $\lambda$ then
$Z_\lambda(\gamma)=\emptyset$. If $\gamma$ is loxodromic and
$l<\lambda$ then $Z_\lambda(\gamma)$ is a cylinder about the axis of
$\gamma$. Its radius
$R_\gamma(\lambda)$ is given explicitly by \Equation\label{omega}
R_\gamma(\lambda)=\arcsinh\bigg(\bigg( \frac{\cosh\lambda -\cosh
  l}{\cosh l-\cos\theta}\bigg)^{1/2}\bigg), \EndEquation where
$\theta$ denotes the twist angle of $\gamma$. This follows from the
formula 
\Equation\label{fred}
\cosh\dist(z,\gamma\cdot z)=\cosh l+(\sinh^2 r)(\cosh
l-\cos\theta)
\EndEquation
which holds for any point $z\in\HH^3$ at a distance $r$ from the axis of a
loxodromic isometry $\gamma$ with translation length $l$ and twist
angle $\theta$.

It follows from the definitions that
$$\tMthin(\lambda )=\bigcup_{1\ne\gamma\in\Gamma} Z_\lambda(\gamma).$$

It follows from our description of the sets $Z_\lambda(\gamma)$ that
the set $Z_\lambda(C)\doteq\bigcup_{1\ne\gamma\in C}Z_\lambda(\gamma)$
is empty for every $C\in\calc(\Gamma)-\calc_\lambda(\Gamma)$; is a
cylinder about $A_C$ for every loxodromic $C\in\calc_\lambda(\Gamma)$;
and is a horoball based at $\omega_C$ for every parabolic $C$.
Furthermore, we have

\Equation\label{union of cylinders}\tMthin(\lambda)=
\bigcup_{C\in\calc_\lambda(\Gamma)} Z_{\lambda}(C),
\EndEquation

We define a {\it Margulis number} for $M$ to be a positive real
number $\lambda$ such that for any two subgroups $C,C'\in\calc_\lambda(\Gamma)$
we have $Z_\lambda(C)\cap Z_\lambda(C')=\emptyset$.  If $\lambda$
is a Margulis number for $M$, the components of $\Mthin(\lambda)$ are
tubes and standard cusp neighborhoods.

Since for any $\lambda>0$ the set $\tMthin(\lambda)$ is a union of sets
of the form $Z_\lambda(\gamma)$, each of which is precisely invariant
under a maximal abelian subgroup, $\lambda$ is a Margulis number if
and only if the fundamental group of each component of
$\Mthin(\lambda)$ has abelian image in $\pi_1(M)$.
\EndNumber

\Number\label{karl} Note also that if $\lambda$ is a Margulis number
then $\tMthick(\lambda)$ is connected, since the union in \ref{union
  of cylinders}
is disjoint. In particular it follows that $\Mthick(\lambda)$ is
connected.   \EndNumber

\Number\label{log 3 works} It follows from \cite[Theorem
10.3]{paradoxical}, together with the main result of \cite{agol} or
\cite{cg}, that if $M$ is a complete, orientable hyperbolic
$3$-manifold without cusps and $\pi_1(M)$ is $2$-free, then $\log3$ is
Margulis number for $M$.
\EndNumber

\section{Cup products and homology of covering spaces}\label{cup product section}

The main result of this section is Proposition \ref{homology stuff},
which is a variation on \cite[Lemma \twoareminusthree]{last}.

\Number 
As in \cite{last}, if $X$
is a topological space, we will set $\mfldrk X = \newrk H_1(X;\ZZ_2)$.
\EndNumber

\Number
As in \cite[\exactsequence]{last}, if $N$ is a normal subgroup of a
group $G$, we shall denote by $G\#N$ the subgroup of $G$ generated by
all elements of the form $gag^{-1}a^{-1}b^2$ with $g \in G$ and $a,b
\in N$.

Also recall from \cite[\exactsequence]{last} the five term exact
sequence due to Stallings \cite{St}.  If $N$ is a normal
subgroup of a group $G$ and if $Q=G/N$, then there is an exact
sequence
\Equation\label{exact sequence} H_2(G;\ZZ_2)\longrightarrow
H_2(Q;\ZZ_2)\longrightarrow N/(G\#N)\longrightarrow
H_1(G;\ZZ_2)\longrightarrow H_1(Q;\ZZ_2)\longrightarrow0,
\EndEquation
in which the maps $H_1(G;\ZZ_2)\longrightarrow H_1(Q;\ZZ_2)$ and
$H_2(G;\ZZ_2)\longrightarrow H_2(Q;\ZZ_2)$ are induced by the natural
homomorphism from $G$ to $Q=G/N$.
\EndNumber

\Number\label{gee sub em} 
As in \cite[\geesubem]{last}, for any group $\Gamma$, we define
subgroups $\Gamma_{d}$ of $\Gamma$ recursively for ${d}\ge0$, by
setting $\Gamma_0=\Gamma$ and $\Gamma_{{d}+1}=\Gamma\#\Gamma_{d}$. We
regard $\Gamma_{d}/\Gamma_{{d}+1}$ as a $\ZZ_2$-vector space.
\EndNumber

The following result is a variation on \cite[Lemma 1.3]{sw}.

\Lemma\label{trivial cup}Suppose that $M$ is a closed,
aspherical $3$-manifold.  Let $t$ denote the rank of the subspace
 of $H^2(M;\ZZ_2)$  spanned by
the image of the
cup product $H^1(M;\ZZ_2)\times H^1(M;\ZZ_2)\to H^2(M;\ZZ_2)$. Set $\Gamma=\pi_1(M)$ and $r=\mfldrk M$.  Then the
$\ZZ_2$-vector space $\Gamma_1/\Gamma_2$ has dimension exactly
$(r(r+1)/2)-t$. 
\EndLemma

\Proof
The group $Q=\Gamma/\Gamma_1$ is an elementary $2$-group of rank
$r$. If we regard $Q$ as a direct product $C_1\times\cdots\times C_r$,
where each $C_i$ is cyclic of order $2$, then the K\"unneth theorem for
cohomology with field coefficients gives an identification 
$$H^2(Q;\ZZ_2)=\bigoplus_{D\in I}V_D,$$
where $I$ is the set of all
$r$-tuples of non-negative integers whose sum is $2$, and
$$V_D=H^{d_1}(C_1;\ZZ_2)\otimes\cdots\otimes H^{d_r}(C_r;\ZZ_2)$$
for
every $D=(d_1,\ldots,d_r)\in I$. We have
$$I=\{D_{hk}:1\le h<k\le r\}\cup\{E_{h}:1\le h\le r\},$$
where the $r$-tuples $D_{hk}$ and $E_h$ are defined by
$D_{hk}=(\delta_{ih}+\delta_{ik})_{1\le i\le r}$ and
$D_{h}=(\delta_{ih})_{1\le i\le r}$. In particular we have
$\dim H^2(Q;\ZZ_2)=|I|=r(r+1)/2$.

Note that when $1\le h<k\le r$, the generator of $V_{D_{hk}}$ is the
cup product of the generators of $H^1(C_h;\ZZ_2)$ and $H^1(C_k;\ZZ_2)$;
and that when $1\le h\le r$, the generator of $V_{D_{h}}$ is the cup
product of the generator of $H^1(C_h;\ZZ_2)$ with itself. (We remark
that the second statement would not hold if we were working over
$\ZZ_p$ for an odd prime $p$, since the cup product of a generator of
$H^1(C_h;\ZZ_p)$ with itself would vanish.)  Hence if
$c:H^1(Q;\ZZ_2)\otimes H^1(Q;\ZZ_2)\to H^2(Q;\ZZ_2)$ is the map defined
by cup product, the image of $c$ contains all the $V_{D_{hk}}$ and all
the $V_{D_h}$, so that $c$ is surjective.

There is a natural commutative diagram
$$\xymatrix{
H^1(Q;\ZZ_2)\otimes H^1(Q;\ZZ_2) \ar@{->}[r]^{\hskip 30pt c} \ar@{->}^{j}[d] &
 H^2(Q;\ZZ_2)\ar@{->}^{\beta^*}[d]\\
H^1(\Gamma;\ZZ_2)\otimes H^1(\Gamma;\ZZ_2) \ar@{->}^{\hskip 30pt \bar c}[r] 
& H^2(\Gamma;\ZZ_2)}
$$
in which the map $\bar c$ is defined by the cup product in the
$\ZZ_2$-cohomology of $\Gamma$,  and $j$ is the 
isomorphism induced by the quotient map $\Gamma\to Q=\Gamma/\Gamma_1$. Since $M$ is aspherical, the $\ZZ_2$-cohomology ring of
$\Gamma$ is isomorphic to that of $M$.  The definition of $t$
therefore implies that the image of $\bar c$ is a $\ZZ_2$-vector space
of dimension $t$. Since we have shown that $c$ is surjective, and
since $j$ is an isomorphism, it follows that the linear map $\beta^*$,
which is induced by the quotient map $\Gamma\to Q$, has rank $t$.
Since $\ZZ_2$ is a field, it follows that the homomorphism
$\beta:H_2(\Gamma;\ZZ_2)\to H_2(Q;\ZZ_2)$ induced by the quotient map
also has rank $t$. 

Applying (\ref{exact sequence}) with $G=\Gamma$ and $N=\Gamma_1$, so
that $Q=\Gamma/\Gamma_1$ is the rank-$r$ elementary $2$-group defined
above and
$N/(G\#N)=\Gamma_1/\Gamma_2$, we obtain an exact sequence
$$
\xymatrix{
H_2(\Gamma;\ZZ_2)\ar@{->}[r]^\beta &
H_2(Q;\ZZ_2)\ar@{->}[r] &
\Gamma_1/\Gamma_2\ar@{->}[r] &
H_1(\Gamma;\ZZ_2)\ar@{->}[r]^\alpha &
\Gamma/\Gamma_1\ar@{->}[r] &
0} .
$$
But we have seen that  $\beta$ has rank $t$, and it is clear that $\alpha$, which
is induced by the quotient map $\Gamma\to Q=\Gamma/\Gamma_1$, is an isomorphism. Hence
by exactness  we have
$$\newrk(\Gamma_1/\Gamma_2)=\newrk(H_2(Q;\ZZ_2))-t=\newrk(H^2(Q;\ZZ_2))-t
=r(r+1)/2-t.$$
\EndProof

\Proposition\label{homology stuff} Suppose that $M$ is a closed,
aspherical $3$-manifold.  Let $t$ denote the rank of the subspace
 of $H^2(M;\ZZ_2)$  spanned by
the image of the
cup product $H^1(M;\ZZ_2)\times H^1(M;\ZZ_2)\to H^2(M;\ZZ_2)$. Set
$r=\mfldrk M$. Then for any integer $m\ge0$ and any regular
covering $\tM$ of $M$ with covering group $(\ZZ_2)^m$, we have
$$\mfldrk\tM\ge(m+1)r-\frac{m(m+1)}2-t.$$
\EndProposition

\Proof
We set $\Gamma=\pi_1(M)$. According to Lemma \ref{trivial cup},
the $\ZZ_2$-vector space $\Gamma_1/\Gamma_2$ has dimension $r(r+1)/2-t$.

Let $N$ denote the subgroup of $\Gamma$ corresponding to the
regular covering space $\widetilde M$. We have
$\Gamma/N\cong(\ZZ_2)^m$. Hence may write $N = E\Gamma_1$ for some
$(r-m)$-generator subgroup $E$ of $\Gamma$.  It now follows from
\cite[Lemma 1.4]{sw} that
\begin{eqnarray*}
\mfldrk\widetilde M &=&\dim H_1(E\Gamma_1;\ZZ_2)\\
&\ge&\dim (\Gamma_1/\Gamma_2) - \frac{(r-m)(r-m-1)}{2}\\
&=&\frac{r(r+1)}{2} -  t  - \frac{(r-m)(r-m-1)}{2}\\
&=& (m+1)r-\frac{m(m+1)}2-t.
\end{eqnarray*}
\EndProof

\section{Dehn filling and \texorpdfstring{$k$}{k}-freeness}
\label{k-free filling section}

This section contains the proofs of Theorem \ref{dewey} and
Proposition \ref{k-free filling}, which were stated in the
Introduction.

\Proposition\label{same difference}Suppose that $N$ is a compact
$3$-manifold whose boundary is a torus, that $F$ is a compact
orientable $2$-manifold with $\chi(F)<0$, and that $f:F\to N$ is a map
such that $f(\partial F)\subset \partial N$, and such that
$f_\sharp:\pi_1(F)\to\pi_1(N)$ is injective. Then
$f_\sharp:\pi_1(F,\partial F)\to\pi_1(N,\partial N)$ is also
injective.
\EndProposition

\Proof Suppose, to the contrary, that there is a path $a$ in $F$ which
is not homotopic rel endpoints into $\partial F$, but that the path
$f(a)$ is homotopic rel boundary into $\partial N$.  The path $a$
determines a map $h:G\to F$ where $G$ is a certain connected graph
with three edges and two vertices.  If $a$ has both endpoints in the
same component $b$ of $\partial F$ then $G$ is a ``theta'' graph; the
restriction of $h$ to one edge is the path $a$ and the restriction to
the union of the other two edges is a homeomorphism onto $b$.  If $a$
has endpoints contained in two components $b_1$ and $b_2$ of $F$ then
$G$ is an ``eyeglass'' graph; the restriction of $h$ to the separating
edge is the path $a$ and the restriction to the union of the other two
edges is a homeomorphism onto $b_1\cup b_2$.  Since $F$ has genus $g >
0$, the map $h$ induces an injection on $\pi_1$.  But the composition
$f\circ h$ is homotopic into $\partial N$, which is a torus, and hence
the induced map on $\pi_1$ is not injective.  It follows that
$f_\sharp:\pi_1(F)\to\pi_1(N)$ is not injective, contradicting the
hypothesis.
\EndProof

\Proof[Proof of Theorem \ref{dewey}] Set $M=N(\alpha)$. It follows
from the hypothesis that $M$ has a connected covering space $p:\tM\to
M$ such that $\pi_1(\tM)$ is isomorphic to the fundamental group
of a closed genus-$g$ orientable surface $\Sigma_g$. The compact core
\cite{compactcore} of $\tM$ is a compact, orientable, irreducible
$3$-manifold $K$ with $\pi_1(K)\cong\pi_1(\Sigma_g)$. Hence $K$ is
homeomorphic to $\Sigma_g\times[-1,1]$ by \cite[Theorem 10.6]{Hempel}. In
particular, there is a genus-$g$ surface $S\subset\tM$ such that

\begin{enumerate}[(1)]
\item\label{huey} the inclusion induces an isomorphism from $\pi_1(S)$ to
  $\pi_1(\tM)$.
\end{enumerate}

We may write $N(\alpha)=N\cup T$, where $T$ is a solid torus 
with $T\cap M=\partial T=\partial M$, and the meridian disks of $T$
represent the slope $\alpha$. Set $\tT=p^{-1}(T)$ and
$\tN=p^{-1}(N)$. We may choose $S$  so that

\begin{enumerate}[(1)]
\setcounter{enumi}{1}
\item\label{louie} each component of $S\cap\tT$ is a meridian disk in $\tT$.
\end{enumerate}

We let $m$ denote the number of components of $S\cap\tT$, and we
suppose $S$ to be chosen, among all genus-$g$ closed surfaces in $\tM$
for which (\ref{huey}) and (\ref{louie}) hold, so that $m$ is as small
as possible. If $m=0$ then $h(S)\subset \tN$, and so the fundamental
group of the connected covering space $\tN$ has a subgroup isomorphic
to a genus-$g$ surface group. It follows that $\pi_1(N)$ also has such a
subgroup in this case.

Now suppose that $m>0$. Set $F=S\cap \tN$, so that $F$ is a compact
orientable surface of genus $g$, properly embedded in $\tN $, with
$\partial F\ne\emptyset$. If the inclusion homomorphism
$\pi_1(F)\to\pi_1(\tN )$ has a non-trivial kernel, it follows from
Dehn's lemma and the loop theorem that there is a disk $D_1\subset \tN
$ such that $D_1\cap F =\partial D_1$, and such that $\partial D_1$
does not bound a disk in $F$. However, it follows from (\ref{huey})
that $\partial D_1$ does bound a disk $D\subset S$, which must contain
at least one component of $S\cap\tT$. If we set
$S_1=\overline{S-D}\cup D_1$, then (\ref{huey}) and (\ref{louie})
still hold when $S$ is replaced by $S_1$, but $S_1\cap\tT$ has at most
$m-1$ components. This contradiction to the minimality of $m$ proves
that $\pi_1(F)\to\pi_1(\tN)$ is injective.

Hence if we set $f=p|F:F\to N$, then $f_\sharp:\pi_1(F)\to\pi_1(N)$ is
injective.  But since $g>0$ we have $\chi(F)<0$, and it therefore
follows from Proposition \ref{same difference} that
$f_\sharp:\pi_1(F,\partial F)\to\pi_1(N,\partial N)$ is also
injective. Thus condition (i) of Definition \ref{singular boundary
slope def} holds.  Since $F=S\cap\tN$, condition (ii) of
\ref{singular boundary slope def} follows from (\ref{louie}). Hence
$\alpha$ is a genus-$g$ singular boundary slope in this case.
\EndProof

\Remark While Theorem \ref{dewey} is similar to \cite[Proposition
1.1]{oertel}, neither of the two results contains the other.  Oertel's
result assumes a weaker hypothesis on the closed surface $F$, since it
only requires that essential {\it simple} closed curves and arcs on
$S$ have homotopically non-trivial images.  On the other hand, the
conclusion of \cite[Proposition 1.1]{oertel} is also weaker since it
does not imply the injectivity of the induced map on relative
fundamental groups.  \EndRemark 

\Proposition\label{variant} Let $M$ be a simple (\ref{simple def}),
closed, orientable $3$-manifold, and let $k\ge2$ be an
integer. Suppose that for some prime $p$ we have
$\dim_{\ZZ_p}H_1(M;\ZZ_p)\ge k+2$, and that $\pi_1(M)$ has no subgroup
isomorphic to a genus-$g$ surface group for any integer $g$ with
$k/2<g<k$. Then $\pi_1(M)$ is $k$-free.  \EndProposition

\Proof This is the same result as \cite[Proposition 1.8, case
(b)]{sw}, except that in the latter result, instead of assuming only
that $\pi_1(M)$ has no subgroup isomorphic to a genus-$g$ surface
group for any integer $g$ with $k/2<g<k$, one makes the superficially
stronger assumption that $\pi_1(M)$ has no subgroup isomorphic to a
genus-$g$ surface group for any integer $g$ with $0<g<k$. However, if
$\pi_1(M)$ does have a subgroup $G$ isomorphic to a genus-$g$ surface
group for some integer $g$ with $0<g<k$, then we must have $g>1$ since
$M$ is simple; and $G$ has a subgroup $H$ of index
$$d=\bigg[\frac{k-2}{g-1}\bigg]\ge1,$$ where $[\cdot]$ is the greatest
integer function. It follows that $H\le\pi_1(M)$ is a genus-$g'$
surface group, where $g'=d(g-1)+1$ satisfies $k/2<g'<k$. This
contradicts the hypothesis of Proposition \ref{variant}.  \EndProof

\Proposition\label{cup product variant}Let $M$ be a simple,
closed, orientable $3$-manifold, let $k\ge2$ be an integer, and
suppose that $\pi_1(M)$ has no subgroup isomorphic to a genus-$g$
surface group for any integer $g$ with $k/2<g<k$. Suppose also that $
\dim_{\ZZ_2}H_1(M;\ZZ_2)\ge k+1$, and the subspace of $H^2(M;\ZZ_2)$
spanned by the image of the cup product $H^1(M;\ZZ_2)\times
H^1(M;\ZZ_2)\to H^2(M;\ZZ_2)$ has  dimension  at most $k-1$.  Then $\pi_1(M)$
is $k$-free.  \EndProposition

\Proof
Set $r= \dim_{\ZZ_p}H_1(M;\ZZ_p)$ and suppose that $\Phi$ is a
subgroup of $\pi_1(M)$ whose rank is at most $k$. Since $r\ge k+1$,
some index-$2$ subgroup $\tGamma$ of $\pi_1(M)$ contains $\Phi$. Let
$\tM$ denote the covering space of $M$ defined by $\tGamma$. It
follows from the case $m=1$ of Proposition \ref{homology stuff} that
$$\dim_{\ZZ_2}H_1(\tM;\ZZ_2)\ge2r-1-t,$$
where $t$ denotes the rank of the subspace of $H^2(M;\ZZ_2)$ spanned
by the image of the cup product $H^1(M;\ZZ_2)\times H^1(M;\ZZ_2)\to
H^2(M;\ZZ_2)$. Since by hypothesis we have $t\le k-1$, it follows that
$$\dim_{\ZZ_2}H_1(\tM;\ZZ_2)\ge2r-k\ge k+2.$$
It therefore follows from Proposition \ref{variant} that
$\tGamma\cong\pi_1(\tM)$ is $k$-free. Hence $\Phi$ is a free group.
\EndProof

\Proof[Proof of Proposition \ref{k-free filling}]
If $\alpha$ is any slope  on $\partial N$, let us denote by
$N(\alpha)$ the closed manifold obtained from $N$ by the Dehn filling
corresponding to $\alpha$.

Let us fix a prime $p$ such that either
\Alternatives
\item $\dim_{\ZZ_p}H_1(N;\ZZ_p)\ge k+2$ for some prime $p$, or
\item $p=2$, the ${\ZZ_2}$-vector space $H_1(N;\ZZ_2)$ has dimension
  at least $ k+1$, and the subspace of $H^2(N;\ZZ_2)$ spanned by the
  image of the cup product $H^1(N;\ZZ_2)\times H^1(N;\ZZ_2)\to
  H^2(N;\ZZ_2)$ has rank at most $k-2$.
\EndAlternatives

Let us fix a basis $\{\lambda,\mu\}$ of $H_1(\partial N;\ZZ)$ such
that $\lambda$ lies in the kernel of the natural homomorphism
$H_1(\partial N;\ZZ)\to H_1(N;\ZZ_p)$. For every integer $n$, let
$\alpha_n$ denote the slope defined by the primitive class
$\lambda+pn\mu\in H_1(\partial N;\ZZ)$.  If $i:N\to N(\alpha_n)$
denotes the inclusion homomorphism, then $i_*:H_1(N;\ZZ_p)\to
H_1(N(\alpha_n);\ZZ_p)$ and $i^*:H^1(N(\alpha_n);\ZZ_p)\to
H^1(N;\ZZ_p)$ are isomorphisms, whereas $i^*:H^2(N(\alpha_n);\ZZ_p)\to
H^2(N;\ZZ_p)$ is surjective and has a $1$-dimensional kernel.

It follows from \cite[Theorem 1]{hass-wang-zhou}, for any given
integer $g>0$ there are only finitely many genus-$g$ singular boundary
slopes in $\partial N$. Hence there is an integer $n_0\ge0$ such that
for every $n$ with $|n|\ge n_0$ and every $g$ with $k/2<g<k$, the
slope $\alpha_n$ fails to be a genus-$g$ singular boundary
slope. Since by hypothesis $\pi_1(N)$ has no subgroup isomorphic to a
genus-$g$ surface group for any integer $g$ with $k/2<g<k$, it follows
from Theorem \ref{dewey} that when $|n|\ge n_0$ and $k/2<g<k$, the
group $\pi(N(\alpha))$ has no subgroup isomorphic to a genus-$g$
surface group.

If (a) holds, then
$$\dim_{\ZZ_p}H_1(N(\alpha_n);\ZZ_p)=\dim_{\ZZ_p}H_1(N;\ZZ_p)\ge k+2$$
for every $n$.  It therefore follows from Proposition \ref{variant}
(with $M=N(\alpha_n)$) that $\pi_1(N(\alpha_n))$ is $k$-free whenever
$|n|\ge n_0$.

Now suppose that (b) holds. Then 
$$\dim_{\ZZ_2}H_1( N(\alpha_n);\ZZ_2)=\dim_{\ZZ_2}H_1( N;\ZZ_2)\ge
k+1$$ for every $n$. Consider the commutative diagram
$$\xymatrix{
H^1(N(\alpha_n);\ZZ_2)\otimes H^1(N(\alpha_n);\ZZ_2)
\ar@{->}[r]^{\hskip 30pt {c_n}} \ar@{->}^{i^*\otimes i^*}[d] &
 H^2(N(\alpha_n);\ZZ_2)\ar@{->}^{i^*}[d]\\
H^1(N;\ZZ_2)\otimes H^1(N;\ZZ_2) \ar@{->}^{\hskip 30pt  c}[r] 
& H^2(N;\ZZ_2)}
$$
in which the homomorphisms $c$ and $c_n$ are defined by cup product. 
Since ${i^*\otimes i^*}$ is an isomorphism, we have 
$$c(H^1(N;\ZZ_2)\otimes H^1(N;\ZZ_2))={i^*}\circ
c_n(H^1(N(\alpha_n);\ZZ_2)\otimes H^1(N(\alpha_n);\ZZ_2)).$$
Since $i^*$  is surjective and
has a $1$-dimensional kernel, it follows that 
$$\dim_{\ZZ_2}c(H^1(N;\ZZ_2)\otimes H^1(N;\ZZ_2))=\dim_{\ZZ_2}
c_n(H^1(N(\alpha_n);\ZZ_2)\otimes H^1(N(\alpha_n);\ZZ_2))-1.$$
Condition (b) gives 
$$\dim_{\ZZ_2}c(H^1(N;\ZZ_2)\otimes H^1(N;\ZZ_2))\le k-2,$$ 
and hence
$$c_n(H^1(N(\alpha_n);\ZZ_2)\otimes H^1(N(\alpha_n);\ZZ_2))\le k-1.$$
It now follows from Proposition \ref{cup product variant} (with
$M=N(\alpha_n)$) that
$\pi_1(N(\alpha_n))$ is $k$-free whenever $|n|\ge n_0$. 
\EndProof

\Remark\label{hatcher} In the proof given above, we quoted
Theorem 1 of Hass, Wang and Zhou's paper \cite{hass-wang-zhou} for
the fact that there are only finitely many genus-$g$ singular boundary
slopes in $\partial N$ for any given $g>0$. Restricting the genus is
essential here: there exist examples \cite{baker} of compact, orientable,
irreducible $3$-manifolds whose boundary is a single torus, in which
there are infinitely many singular boundary slopes.

The result proved by Hatcher in \cite{hatcher} implies that if one
considers only embedded surfaces, rather than singular ones, a
finiteness result holds without restricting genus.  A slope $\beta$ in
$\partial N$, where $N$ is a compact, orientable, irreducible
$3$-manifold whose boundary is a torus, is called a {\it boundary
slope} if there is a connected orientable surface $F$ with non-empty
boundary, properly embedded in $N$ and not boundary-parallel, such
that the inclusion homomorphism $\pi_1(F)\to\pi_1(N)$ is injective and
all components of $\partial F$ represent the slope $\beta$. The main
theorem of \cite{hatcher} implies that there are only finitely many
boundary slopes in $\partial N$.
\EndRemark

\section{\texorpdfstring{$k$}{k}-free groups and generating sets}
\label{algebra section}

This brief section provides the algebra needed for the material in the
following sections.

We shall denote by $|X|$ the cardinality of a set $X$. If $S$ is a
subset of a group $\Gamma$ we shall denote by $\langle S\rangle$ the
subgroup of $\Gamma$ generated by $S$.

We shall say that a subset $S$ of a group $\Gamma$ is {\it
  independent} (or that the elements of $S$ are independent) if
$\langle S\rangle$ is a free group with basis $S$.

\Number\label{homer} It is a basic fact in the theory of free groups
\cite[Vol. 2, p. 59]{kurosh} that a finite set $S\subset\Gamma$ is
independent if and only if $\langle S\rangle$ is free of rank $|S|$.
\EndNumber 

\Proposition\label{algebra} Let $k$ be a positive integer. Suppose
that $\Gamma$ is a finitely generated group which is $k$-free but not
free, that $X$ is a generating set for $\Gamma$, and that $T$ is a
finite independent subset of $X$ with $|T|\le k$. Then there is
an independent subset $S$ of  $X$  such that $T\subset S$
and $|S|=k$.  \EndProposition

\Proof
Set $m=|T|$. If $m=k$ we may take $S=T$. If $m<k$ we shall show that
there is an independent subset $T'$ of $\Gamma$ such that $T\subset
T'\subset X$ and $|T'|=m+1$. The result will then follow at once by
induction.

Let us write $X=\{x_1,\ldots,x_n\}$, where $x_1,\ldots,x_n$ are
distinct and are indexed in such a way that $T=\{x_1,\ldots,x_m\}$.
For each $j$ with $0\le j\le n$, let $A_j$ denote the subgroup of
$\Gamma$ generated by $x_1,\ldots, x_j$. Since $A_n=\Gamma$ is $k$-free but
not free, it must have rank $>k$. On the other hand, since $T$ is
independent, $A_m=\langle T\rangle$ is a free group of rank $m<k$.
Hence there is an index $p$ with $m\le p<n$ such that $A_p$ has rank
$\le m$ and $A_{p+1}$ has rank $>m$. Set $\xi=x_{p+1}$. Since
$A_{p+1}=\langle A_p,\xi\rangle$ we must have $\rank A_{p+1}\le1+\rank
A_p$. Hence $\rank A_{p}=m$ and $\rank A_{p+1}=m+1$.

Since $\Gamma$ is $k$-free and $m<k$, the subgroups $A_p$ and
$A_{p+1}$ are free. Let $y_1,\ldots,y_m$ be a basis for $A_p$. Then
the rank-$(m+1)$ free group $A_{p+1}$ is generated by the elements
$y_1,\ldots,y_m,\xi$, which must therefore form a basis for $A_{p+1}$
(see \ref{homer}).  This shows that $A_{p+1}$ is a free product
$A_p*\langle\xi\rangle$, and that the factor $\langle\xi\rangle$ is
infinite cyclic. Hence the subgroup $\langle A_m,\xi\rangle$ of
$A_{p+1}$ is a free product $A_m*\langle\xi\rangle$. Since
$T=\{x_1,\ldots,x_m\}$ freely generates $A_m$ it follows that the set
$T'=T\cup\{\xi\}$ has cardinality $m+1$ and freely generates
$A_{p+1}$. In particular, $T'$ is independent.
\EndProof

\section{Distant points and \texorpdfstring{$k$}{k}-freeness}
\label{distant section}

The two main results of this section, Proposition \ref{big
radius} and Proposition \ref{popeye}, were discussed in the
Introduction.

\Proposition\label{big radius}Let $k$ and $m$ be integers $k\ge2$ and
$0\le m\le k$. Suppose that $M=\HH^3/\Gamma$ is a closed, orientable
hyperbolic $3$-manifold such that $\pi_1(M)\cong\Gamma$ is
$k$-free. Let $q:\HH^3\to M$ denote the quotient map. Suppose that
$N\subset M$ is a closed set such that $q^{-1}(N)$ is connected.  Let
$P$ be a point of $M$, let $\xi_1,\ldots,\xi_m$ be independent
elements of $\pi_1(M,P)$ represented by loops $\ell_1,\ldots,\ell_m$
based at $P$, and let $\lambda_j$ denote the length of $\ell_j$.  Then
there is a point $Q\in N$ such that $\rho=\dist_M(P,Q)$ satisfies
$$\frac{k-m}{1+e^{2\rho}}+
\sum_{j=1}^m\frac{1}{1+e^{\lambda_j}}
\le\frac12.$$
\EndProposition

\Proof Set
$$\rho=\max_{Q\in N}\dist_M(P,Q).$$
Fix a point $\widetilde P\in\HH^3$ with $q(\widetilde P)=P$. Then
$N\subset q(B)$, where $B$ denotes the open ball of radius $\rho$
centered at $\widetilde P$. Hence if we set $\widetilde N=q^{-1}(N)$,
we have
$$\widetilde N\subset\Gamma\cdot
B\doteq\bigcup_{\gamma\in\Gamma}\gamma(B).$$

Let $X_0$ denote the set of all elements $\eta\in\Gamma$ such that
$\dist(\eta(\widetilde P), \widetilde P)<2\rho$.  Note that the
triangle inequality gives $\eta\in X_0$ for any $\eta\in\Gamma$ such
that $\eta(B)\cap B\ne\emptyset$.

We claim that $X_0$ generates $\Gamma$. To show this, we set
$\Gamma_0=\langle X_0\rangle\le\Gamma$, and choose a system $Y$ of
left coset representatives for $\Gamma_0$ in $\Gamma$. For each $y\in
Y$, we set $\widetilde N_y=\widetilde N\cap y\Gamma_0\cdot B$, where
$$y\Gamma_0\cdot B=\bigcup_{\gamma\in\Gamma}y\gamma(B).$$
Then
$$\widetilde N=\bigcup_{y\in Y}\widetilde N_y.$$

If $y$ and $y'$ are distinct elements of $Y$ such that $\widetilde
N_y\cap\widetilde N_{y'}\ne\emptyset$, then for some
$\gamma,\gamma'\in\Gamma_0$ we have $y\gamma(B)\cap
y'\gamma'(B)\ne\emptyset$, so that $B\cap
\gamma^{-1}y^{-1}y'\gamma'(B)\ne\emptyset$. It follows that
$\gamma^{-1}y^{-1}y'\gamma'\in \Gamma_0$ and hence that
$y^{-1}y'\in\Gamma_0$, a contradiction since $y$ and $y'$ represent
distinct left cosets of $\Gamma_0$. This shows that the sets
$\widetilde N_y$ are pairwise disjoint as $y$ ranges over $Y$. Since
these sets are open in the subspace topology of the connected set
$\widetilde N$, we must have $\widetilde N_y=\widetilde N$ for some
$y\in Y$. Since $P\in N$ it follows that $\widetilde N_y$ contains the
$\Gamma$-orbit of $\widetilde P$; and as $\Gamma$ acts freely on
$\HH^3$, this implies that $\Gamma_0=\Gamma$, as claimed.

By hypothesis there exist $m$ independent elements
$\xi_1,\ldots,\xi_m$ of $\pi_1(M,P)$ which are represented by
piecewise smooth loops of respective lengths
$\lambda_1,\ldots,\lambda_m$ based at $P$. If we choose a point
$\widetilde P\in q^{-1}(P)$, we may regard $(\HH^3,\widetilde P)$ as a
based covering space of $(M,P)$, and identify t%he deck transformation
group $\Gamma$ with $\pi_1(M,P)$.  Then for $j=1,\ldots,m$ we have
$\dist(\widetilde P,\xi_j(\widetilde P))\le\lambda_j$.

Let us set $X=X_0\cup\{\xi_1,\ldots,\xi_m\}$.  We wish to apply
Proposition \ref{algebra}, taking $k$ and $\Gamma$ to be defined as
above, and taking $T=\{\xi_1,\ldots,\xi_m\}$ and $X=T\cup X_0$.  By
hypothesis the finitely generated group $\Gamma\cong\pi_1(M)$ is
$k$-free. To see that it is not free, we need only observe that
$H_3(\Gamma;\ZZ)$ has rank $1$ because $M$ is orientable, closed and
aspherical.  Since $X_0$ generates $\Gamma$, in particular the set
$X=\{\xi_1,\ldots,\xi_m\}\cup X_0$ generates $\Gamma$.  We have
$T\subset X$ by definition, and $T$ is independent by hypothesis.
Hence Proposition \ref{algebra} gives an independent subset $S$ of
$\Gamma$ such that $T\subset S\subset X$ and $|S|=k$.  We may write
$S=\{\xi_1,\ldots,\xi_{k}\}$, where $\xi_{m+1},\ldots,\xi_k$ are
elements of $X_0$.

Since $S$ is independent, it follows from \cite[Theorem 6.1]{accs},
 together with the main result of \cite{agol} or \cite{cg},  that
$$\sum_{j=1}^{k} \frac{1}{1+e^{d_j}}\leq\frac12,$$
where $d_j=\dist(\widetilde P,\xi_j\cdot \widetilde P)$ for $j=1,\ldots,k$.

Since $\xi_{m+1},\ldots,\xi_k$ belong to $X_0$, we have $d_j\le2\rho$ for
$j=m+1,\ldots,k$. We have seen that $d_j\le\lambda_j$ for $j=1,\ldots,m$.
Hence
$$\frac{k-m}{1+e^{2\rho}}+
\sum_{j=1}^m\frac{1}{1+e^{\lambda_j}}
\le\frac12,$$
and the proposition is proved.
\EndProof   

\Corollary\label{big radius corollary}
Let $k$ and $m$ be integers $k\ge2$ and $0\le m\le k$. Suppose that
$M=\HH^3/\Gamma$ is a closed, orientable hyperbolic $3$-manifold such
that $\pi_1(M)\cong\Gamma$ is $k$-free, and let $\mu$ be a Margulis
number for $M$.  Let $P$ be a point of $M$, let $\xi_1,\ldots,\xi_m$
be independent elements of $\pi_1(M,P)$ represented by loops
$\ell_1,\ldots,\ell_m$ based at $P$, and let $\lambda_j$ denote the
length of $\ell_j$.  Then there is a point $Q\in \Mthick(\mu)$
such that $\rho=\dist_M(P,Q)$ satisfies
$$\frac{k-m}{1+e^{2\rho}}+
\sum_{j=1}^m\frac{1}{1+e^{\lambda_j}}
\le\frac12.$$
\EndCorollary

\Proof Write $M=\HH^3/\Gamma$, and let $q:\HH^3\to M$ denote the
quotient map.  As we pointed out in \ref{karl}, the set
$\tMthick(\mu)=q^{-1}(\Mthick(\mu))$ is connected. Hence the
assertions of the corollary follow from Proposition \ref{big radius}
if we take $N=\Mthick(\mu)$.
\EndProof

\Proposition\label{popeye} Let $k\ge2$ be an integer, let $M$ be a
closed, orientable hyperbolic $3$-manifold, and suppose that
$\pi_1(M)$ is $k$-free. Let $\mu$ be a Margulis number for $M$.  Let
$c$ be a closed geodesic in $M$, let $l$ denote the length of $c$, and
suppose that $l<\mu$.  Let $R$ denote the radius of the maximal
(embedded) tube about $c$, and set
  $$R'=\frac12\arccosh\big(\cosh  (2R)  \cosh\big(\frac l2\big)\big).$$
Then  there exist points $P\in c$ and $Q\in\Mthick(\mu)$ such that 
$\rho=\dist_M(P,Q)$
satisfies
$$\frac{k-2}{1+e^{2\rho}}+\frac{1}{1+e^{l}}+\frac{1}{1+e^{2R'}}
\le\frac12.$$
\EndProposition

\Proof 
Let us choose a component $A$ of $q^{-1}(c)$. Then $A$ is the axis of
a loxodromic isometry $\gamma_0\in\Gamma$, which has translation
length $l$ and generates a maximal cyclic subgroup $C_0$ of $\Gamma$.

Since $R$ is the radius of a maximal tube about $c$, there is an
element $\eta$ of $\Gamma-C_0$ such that the open cylinders $N(R,A)$
and $N(R,\eta\cdot A)$ are disjoint but have intersecting
closures. Hence the minimum distance between the lines $A$ and
$\eta\cdot A$ is $2R$. Let $\tP$ and $\tP'$ denote the respective
points of intersection of $\eta\cdot A$ and $A$ with their common
perpendicular, so that $\dist(\tP,\tP')=2R$. Let $P,P'\in M$ denote
the images of $\tP$ and $\tP'$ under the quotient map $\HH^3\to
M$. Since $\eta^{-1}\cdot \tP \in A$, the $C_0$-orbit of
$\eta^{-1}\cdot \tP$ contains a point $\tP''=\gamma_0^m\eta^{-1}\cdot
\tP$ such that $\dist(\tP ',\tP'')\le l/2$. Setting $d=\dist(\tP,\tP'')$,
we have
$$\cosh d=\cosh(2R)\cosh\dist(\tP',\tP'')\le\cosh(2R)\cosh(l/2)$$
and hence $d\le2R'$.

We regard $(\HH^3,\widetilde P'')$ as a based covering space of
$(M,P)$, and identify the deck transformation group $\Gamma$ with
$\pi_1(M,P)$. Since $\eta\notin C_0$, the elements $\eta
\gamma_0^{-m}$ and $\gamma_0$ of $\pi_1(M,P)$ do not commute, and
since $\Gamma$ is in particular $2$-free, they are independent. Since
$\dist(\tP '',\tP )=d$ and $\dist(\tP '',\gamma_0\cdot \tP '')=l$, we
may apply Corollary ref{big radius corollary} with $m=2$,
$\xi_1=\gamma_0$, $\xi_2=\eta \gamma_0^{-m}$, $\lambda_1=l$ and
$\lambda_2=d$, to obtain a point $Q\in\Mthick(\mu)$ such that
$\rho=\dist_M(P ,Q)$ satisfies
$$\frac12\ge
\frac{k-2}{1+e^{2\rho}}+
\frac{1}{1+e^{l}}+\frac{1}{1+e^{d}}\ge
\frac{k-2}{1+e^{2\rho}}+
\frac{1}{1+e^{l}}+\frac{1}{1+e^{2R'}}.$$
\EndProof

\section{Distant points from maximal cusp neighborhoods}\label{moe}

We begin by reviewing  a version of W. Thurston's Dehn Surgery
Theorem and proving a result, Proposition \ref{tube convergence},
which is a consequence of the proof of Thurston's theorem.  The
proofs of Thurston's result that have appeared vary in the level of
detail as well as in the strength of the statement, but most of them
are based on the same ideas.  For the sake of definiteness we refer
the reader to the proof given by Petronio and Porti in
\cite{petronio-porti} which seems to contain the most complete and
elementary treatment of a version of Thurston's theorem that is
adequate for our purposes.

Suppose that $M$ is a a complete, finite-volume, orientable hyperbolic
manifold with cusps and that $M'$ is a closed manifold obtained by
Dehn filling the compact core of $M$.  Following
\cite{benedetti-petronio}, we shall say that $M'$ is a {\it hyperbolic
  Dehn filling} of $M$ provided that $M'$ admits a hyperbolic
structure in which the core curves of the filling solid tori are 
isotopic to geodesics, which we shall call {\it core geodesics}.

While the general versions of Thurston's Theorem apply to finite-volume
hyperbolic $3$-manifolds with arbitrarily many cusps, we are only
concerned here with the case of one cusp.  We record the following
statement, which follows from any of the various versions of Thurston's
Theorem.

\Theorem[Thurston]\label{Dehn filling for the lumpen proletariat} Let
$M$ be a one-cusped complete, finite-volume, orientable hyperbolic
$3$-manifold.  Then all but finitely many Dehn fillings of the compact
core of $M$ are hyperbolic Dehn fillings of $M$.  \EndTheorem

Here is the result whose proof is extracted from that of Thurston's
theorem. The deepest conclusions are (\ref{mr bluster}) and
(\ref{dilly-dally}). Conclusion (3) was used in the proof of
\cite[Lemma 4.3]{AgolCS}.  While we presume that all of the
conclusions are well known to experts in geometric convergence,
conclusion (4) was new to us.  In the proof we have attempted to
follow the terminology used in \cite{jorgensen-marden} as closely as
possible.

\Proposition\label{tube convergence} Let
$M_\infty$ be a complete, finite-volume, orientable hyperbolic
manifold having exactly one cusp and let $(M_j)_{j\in\NN}$ be a
sequence of distinct hyperbolic Dehn fillings of $M_\infty$.  Let
$T_j$ be a maximal tube about the core geodesic for the Dehn filling
$M_j$ and let $H$ be a maximal standard neighborhood of the cusp in
$M$.  Then 
\Conclusions
\item\label{flubadub} $\vol M_j \to \vol M_\infty$ as $j\to \infty$;
\item\label{clarabelle}the length of the core geodesic of $M_j$ tends
  to $0$ as $j\to\infty$;
\item\label{mr bluster}  $\vol T_j \to \vol H$ as $j\to\infty$; 
\item\label{dilly-dally} if a given positive number $\lambda$ has the
  property that it is a Margulis number for each $M_j$, then $\lambda$
  is also a Margulis number for $M_\infty$.
\EndConclusions 
\EndProposition

\Proof To start, we recall some of the statements proved in
\cite{petronio-porti} in the course of proving Thurston's theorem.  We
let $U$ denote the open unit disk in $\CC$ and let
$\rho_0:\pi_1(M_\infty)\to PSL_2(\CC)$ be a discrete faithful
representation.  By a {\it hyperbolic ideal tetrahedron} we mean the
convex hull in $\HH^3$ of a $4$-tuple of distinct points on the sphere at
infinity of $\HH^3$.  (We call these points vertices, although they
are not contained in the tetrahedron.)  We say that a hyperbolic ideal
tetrahedron is {\it flat} if it is contained in a plane.   If the
vertices of a non-flat hyperbolic ideal tetrahedron are given an
ordering, then the cross-ratio of the $4$-tuple of vertices has
positive imaginary part if and only if the ordering is consistent with
the orientation that the tetrahedron inherits from hyperbolic space.
There is a smooth function, which we shall denote $V(w)$, that is
defined on $\CC-\{0,1\}$ and has the property that if the vertices of a
hyperbolic ideal tetrahedron $\Delta$ are ordered, and if $w$ is the
cross-ratio of the $4$-tuple of vertices, then $V(w) = \vol\Delta$ if
$\Im w \ge 0$ and $V(w) = -\vol\Delta$ if $\Im w \le 0$.

For $j\in \NN\cup\{\infty\}$ we give $M_j$ the orientation that it
inherits as a quotient of $\HH^3$.

Let $\widetilde M_\infty$ and $\widetilde M_j$ denote the universal covers of
$M_\infty$ and $M_j$ respectively.  We will set $\Gamma=\pi_1(M_\infty)$.

The results that we need from \cite{petronio-porti} are summarized below:
\Bullets

\item For some $\pi_1(M)$-equivariant topological ideal triangulation
  $\mathcal T$ of $\widetilde M$ there exists a $\rho_0$-equivariant
  developing map $D_0:\widetilde M\to \HH^3$ that  maps each ideal
  $3$-simplex of $\mathcal T$ to a hyperbolic ideal tetrahedron and
  preserves the orientation of each ideal $3$-simplex which does not
  have a flat image. Moreover $D_0$ is a homeomorphism on the union of
  the interiors of the ideal $3$-simplices with non-flat image. 
\medskip

\item $\rho_0$ lies in a continuous family (with respect to the
  complex topology on the representation variety) of
  representations $\rho_z:\Gamma\to PSL_2(\CC)$, for $z\in U$.  This
  family has the property that, given an arbitrarily small
  neighborhood $O$ of $0$, all but finitely many Dehn fillings on $M$
  are realized as hyperbolic Dehn fillings by hyperbolic manifolds of
  the form $\HH^3/\rho_z(\Gamma)$, for $z\in O$.
\medskip

\item for each $z\in U$ there is a $\rho_z$-equivariant developing map
  $D_z:\widetilde M\to\HH^3$ which sends each tetrahedron in $\mathcal T$
  to a hyperbolic ideal tetrahedron.  Moreover the cross-ratios of the
  image tetrahedra vary continously with $z$.

\EndBullets

According to the second result above we may assume,
after deleting finitely many terms from the sequence $(M_j)$,
that there are complex numbers $z_j \in $U such that
$M_j = \HH^3/\rho_{z_j}(\Gamma)$ for $j\in \NN$.

Let $\Delta_1, \ldots, \Delta_n$ be a system of distinct
representatives for the $\Gamma$ orbits of ideal $3$-simplices in
$\mathcal T$.  Order the vertices of each $\Delta_i$ in a way that is
consistent with the orientation inherited from $M_\infty$.  For
$i\in\{1, \ldots, n\}$ and $z\in U$ let $w_i(z)$ denote the
cross-ratio of the vertices of $D_z(\Delta_i)$, with the ordering
determined by the ordering of the vertices of $\Delta_i$.  The results
above then imply that the function
$$\mathcal V(z) = \sum_{i=1}^n V(w_i(z))$$
has the property that $V(z_j) = \vol M_j$.  Since $\mathcal V$ is
continuous, this implies conclusion (\ref{flubadub}).

The continuity of the family $\rho_z$ implies, in the terminology of
\cite{jorgensen-marden}, that the Kleinian groups $\Gamma_j =
\rho_{z_j}(\Gamma)$ converge algebraically to the Kleinian group
$\Gamma_\infty = \rho_0(\Gamma)$.  In particular, there are surjective
homomorphisms $\phi_j:\Gamma_\infty \to \Gamma_j$ such that each
element $\gamma$ of $\Gamma_\infty$ is the limit of the sequence
($\phi_j(\gamma))$.  For each $j\in \NN\cup\{\infty\}$ we shall let
$p_j:\HH^3\to M_j$ denote the quotient map.

We fix a base point in the boundary of a compact core $N$ of $M$, and
let $\Lambda\le\Gamma$ denote the image of $\pi_1(\partial N)$ under
the inclusion homomorphism.  Set $\Lambda_\infty = \rho_0(\Lambda)$
and $\Lambda_j = \rho_j(\Lambda)$ for $j\in\NN$.  Thus
$\Lambda_\infty$ is a maximal parabolic subgroup of $\Gamma_\infty$
and the groups $\Lambda_j$ for $j\in \NN$ are cyclic groups generated
by a loxodromic isometry.

To prove conclusion (\ref{clarabelle}), for $j\in\NN$ we let
$l_j$ denote the length of the core geodesic of $M_j$ for $j\in\NN$
and let $g_j$ denote a loxodromic isometry of translation length $l_j$
that generates $\Lambda_j$.

Let us fix, arbitrarily, a non-trivial element $\gamma_0$ of
$\Lambda$. For every $j\in\NN$ we have $\rho_{z_j}(\gamma_0)=g_j^{m_j}$
for some $m_j\in\ZZ$. Since the Dehn fillings in the sequence are
distinct, there is at most one index $j$ for which
$\rho_{z_j}(\gamma_0)$ is the identity. In particular we have
$m_j\ne0$ for large $j$, so that $\rho_{z_j}(\gamma_0)$ is loxodromic with
translation length $|m_j|l_j\ge l_j$.

On the other hand, algebraic convergence implies that
$\rho_{z_j}(\gamma_0)$ approaches the parabolic isometry
$\rho_{z_j}(\gamma_0)$ as $j\to\infty$. In particular, $(\trace
\rho_{z_j}(\gamma_0))^2\to4$, so that the translation length of
$\rho_{z_j}(\gamma_0)$ approaches $0$. Hence $l_j\to0$, which is
conclusion (\ref{clarabelle}).

Our proofs of conclusions (\ref{mr bluster}) and (\ref{dilly-dally})
depend on the techniques of \cite{jorgensen-marden}. It follows
from \cite[Proposition 3.8]{jorgensen-marden} that algebraic
convergence implies polyhedral convergence, for the sequence of
non-elementary groups $(\Gamma_j)_{j\in\NN}$ as well as for the
sequence of elementary groups $(\Lambda_j)_{j\in\NN}$ (see
\cite[Remark 3.10]{jorgensen-marden}).  We will describe in some
detail what is meant by polyhedral convergence of  the sequence
 $(\Gamma_j)$, and introduce some notation that will be used later.
An analogous description applies to the sequence 
$(\Lambda_j)$.  Consider a point $\widetilde P\in\HH^3$, chosen from
the dense set of points whose associated Dirichlet domains for the
groups $\Gamma_\infty$ and $\Gamma_j$, $j\in\NN$ are generic in the
sense of \cite{jorgensen-marden}.  For $j\in\NN\cup\{\infty\}$ let us
set $P_j = p_j(\widetilde P)$, and $W_j(r)=\overline{N(r, P_j)}\subset
M_j$ for each $r>0$. We denote by $Z_j$  the Dirichlet
domain for $\Gamma_j$ centered at $\widetilde P$.  The truncated
Dirichlet domain $Z_j'=Z_j\cap \overline{N(r, \widetilde P)}$ (a
polyhedral set having some faces that are planar and some that lie in
a sphere) is a fundamental domain for the action of $\Gamma_j$ on
$p_j^{-1}(W_j(r))$.  To say that $(\Gamma_j)_{j\in\NN}$ converges
polyhedrally to $\Gamma_\infty$ means that for any sufficiently large
positive number $r$, say for $r\ge r_1$, there exists $J_r\in\NN$ such
that for $j>J_r$ the polyhedral sets $Z_j'$ and $Z_\infty'$ have
isomorphic combinatorial structures, with face-pairings that
correspond under the homomorphism $\phi_j$.  Moreover, for $j>J_r$,
there are positive numbers $\epsilon_j$ tending to $0$ and
homeomorphisms $g_j:Z_\infty' \to Z_j'$ which preserve vertices, edges
and faces, commute with the face-pairings, and satisfy $\dist(x,
g_j(x)) < \epsilon_j$.

It follows from \cite[Proposition 4.7]{jorgensen-marden} that when
$r\ge r_1$, the set $\overline{M_\infty - W_\infty(r)}$ is a 
closed subset of  $M_\infty$ that is homeomorphic to $S^1\times
S^1\times[1,\infty)$, and that for $j > J_ r$ the set $M_j - W_j(r)$
is a topological regular neighborhood of the core geodesic of the Dehn
filling.

When $r\ge r_1$ it follows from polyhedral convergence that for $j >
J_r$ there exist bi-Lipschitz homeomorphisms $f_{j,r}:W_\infty(r) \to
W_j(r)$ whose Lipschitz constants tend to $1$ as $j\to\infty$. (In
general we shall write $f_j$ for $f_{j,r}$, the choice of $r$ being
clear from the context.) We may take the homeomorphisms $f_j$ to have
the following property:

{
\sl\quad for any $x\in W_\infty(r)$ there exist $\widetilde
x_\infty\in p_\infty^{-1}(x)$ and $\widetilde x_j \in p_j^{-1}(f_j(x))$
with $\dist(\widetilde x_\infty, \widetilde x_j) < \epsilon_j$.
}

Given a subgroup $G$ of $\Gamma_j$ for some $j\in\NN\cup\{\infty\}$ we
define a function $\mu_G:M_j\to\RR$ by $\mu_G(x) =
\min\{\dist(\widetilde x, \gamma(\widetilde x)\}$ where the minimum is
taken over all elements $\gamma$ of $G$ and all points $\widetilde
x\in p_j^{-1}(x)$.  Note that $\mu_{\Gamma_j}(x)$ is twice the
injectivity radius of $M_j$ at $p_j(x)$.

We now claim:
\begin{claim}\label{lisa}  
For any $r \ge r_1$ and for any point $x\in W_\infty(r)$ we have
$\mu_{\Gamma_j}(f_j(x)) \to \mu_{\Gamma_\infty}(x)$
and $\mu_{\Lambda_j}(f_j(x)) \to \mu_{\Lambda_\infty}(x)$ as $j\to\infty$.
\end{claim}

To prove the first assertion of \ref{lisa}, let $r \ge r_1$ be given,
and let $x$ be a given point in $W_\infty(r)$.  Let $\widetilde x$ and
$\widetilde x_j$ be as above.  Then, for any
$\gamma\in\Gamma_\infty-\{1\}$, we have $\widetilde x_j\to \widetilde
x_\infty$ and $\phi_j(\gamma)\to \gamma$ as $j\to\infty$.  Thus
$\dist(\widetilde x_j,\phi_j(\gamma)(\widetilde x_j)) \to
\dist(\widetilde x_\infty, \gamma(\widetilde x_\infty))$ as
$j\to\infty$.  It follows that $\phi_j(\gamma) \not= 1$ for
sufficiently large $j$.  Hence if we choose $\gamma\in
\Gamma_\infty-\{1\}$ so as to minimize $\dist(\widetilde x_\infty,
\gamma(\widetilde x_\infty))$ then for any $\epsilon > 0$ we have
$$\mu_{\Gamma_\infty}(x) = \dist(\widetilde x_\infty,
\gamma(\widetilde x_\infty)) \ge \dist(\widetilde x_j,
\phi_j(\gamma)(\widetilde x_j)) - \epsilon \ge \mu_{\Gamma_j}(f_j(x))
- \epsilon.$$  It follows that $\mu_{\Gamma_\infty}(x) \ge
\limsup \mu_{\Gamma_j}(f_j(x))$.

Next we shall show that $\mu_{\Gamma_\infty}(x) \le \liminf
\mu_{\Gamma_j}(f_j(x))$. It will simplify the notation if we assume,
as we may by passing to a subsequence, that the sequence
$\big(\mu_{\Gamma_j}(f_j(x))\big)$ converges.  For each $j\in\NN$ let
us choose $\gamma_j\in\Gamma_\infty$ such that $\mu_{\Gamma_j}(f_j(x))
= \dist(\widetilde x_j,\phi_j(\gamma_j)(\widetilde x_j))$.  For
sufficiently large $j$ we have $\mu_{\Gamma_j}(f_j(x)) < 1 +
\mu_{\Gamma_\infty}(x)$.  In particular, this shows that the set
$\{\dist(\widetilde x_j,\phi_j(\gamma_j)(\widetilde x_j)) \,|\,
j\in\NN\}$ is bounded and hence that the set $\{\phi(\gamma_j) \,|\,
j\in\NN\}$ is a bounded subset of $PSL(2,\CC)$.  We now appeal to
\cite[Proposition 3.10]{jorgensen-marden} which states that a sequence
of non-elementary Kleinian groups converges polyhedrally if and only
if it converges geometrically, and that the geometric and polyhedral
limits are the same whenever they exist. Thus $\Gamma_\infty$ is the
geometric limit of $(\Gamma_n)$.  Now the definition of geometric
convergence (\cite[Definition 3.2]{jorgensen-marden}) implies in
particular that if we are given any subsequence $\Gamma_{j_k}$ of
$(\Gamma_j)$ and elements $g_{j_k}\in\Gamma_{j_k}$ such that
$g_{j_k}\to h\in PSL(2,\CC)$, then $h$ must be an element of
$\Gamma_\infty$.  Since $\{\phi(\gamma_j) \,|\, j\in\NN\}$ is bounded
we may choose a subsequence $(\gamma_{j_k})$ of $(\gamma_j)$ such that
$\phi(\gamma_{j_k}) \to h\in\Gamma_\infty$.  We then have
$$\mu_{\Gamma_{j_k}}(f_{j_k}(x)) = \dist(\widetilde
x_{j_k},\phi_{j_k}(\gamma_{j_k})(\widetilde x_{j_k})) \to
\dist(\widetilde x, \gamma(x)) \ge \mu_{\Gamma_\infty}(x).$$ Thus,
under our assumption that $(\mu_{\Gamma_j}(f_j(x))$ converges, we have
shown that $\mu_{\Gamma_\infty}(x)\le
\lim_{j\to\infty}\mu_{\Gamma_j}(f_j(x))$.  Dropping this assumption, we
conclude in the general case that $\mu_{\Gamma_\infty}(x)\le
\liminf_{j\to\infty}\mu_{\Gamma_j}(f_j(x))$, as required to complete the
proof of  the first assertion of \ref{lisa}.

The proof of the second assertion of \ref{lisa} is almost the same as
the proof of the first assertion, with the sequence $(\Gamma_j)$
replaced by the sequence $(\Lambda_j)$. The only difference occurs in
the step where we applied \cite[Proposition 3.10]{jorgensen-marden},
which applies only to a sequence of non-elementary Kleinian
groups. For the proof of the second assertion, we appeal instead to
\cite[Remark 3.10]{jorgensen-marden}, which gives the same conclusion
for the special  sequence $(\Lambda_j)$ of elementary Kleinian
groups. Thus \ref{lisa} is proved.

For the rest of the argument, to unify the notation, we set
$T_\infty=H\subset M_\infty$.  We next claim:

\begin{claim}\label{ricky}
For each $j\in\NN\cup\{\infty\}$, the function 
$\mu_{\Lambda_j}$ takes a constant value $I_j$ on the frontier of 
$T_j$. On the open set $T_j$, $\mu_{\Lambda_j}$ is a
monotone decreasing function of distance from the frontier. 
\end{claim}

To prove \ref{ricky} when $j\in\NN$ we may apply the formula
(\ref{fred}), which implies that for $z\in p_j^{-1}(T_j)$ we have
\Equation\label{ethel} \cosh\mu_{\Lambda_j}(z)=\min_{n\in\NN}(\cosh
nl+(\sinh^2 r)(\cosh nl-\cos n\theta)),
\EndEquation 
where $l$ and $\theta$ denote the twist angle of a generator of
$\Lambda_j$ and $r$ is the distance from $z$ to the axis of
$\Lambda_j$.  Since the factor $\cosh nl-\cos n\theta$ increases
monotonically with $n$ for large $n$, the minimum over $\NN$ in
(\ref{ethel}) may be replaced by a minimum over a finite subset of
$\NN$. Thus $\mu_{\Lambda_j}(z)$ is the minimum of a finite collection
of smooth, monotone increasing functions of $r$, and is therefore
monotone increasing and piecewise smooth. This proves \ref{ricky} for
$j\in\NN$; for $j=\infty$ it is similar but easier.

It follows from \ref{ricky} that with the exception of the core curve
of $T_j$ for $j\in\NN$, the level sets of the $\mu_{\Lambda_j}$ for
$j\in\NN\cup\{\infty\}$ are smooth embedded tori.  For any $r\ge r_1$ and any
$j\in\NN\cup\{\infty\}$ we may characterize the submanifold $T'_j =
T_j\cap W_j(r)$ as the set $\{x\in W_j(r)\,|\, \mu_{\Lambda_j}(x) <
I_j\}$. The frontier of $T'_j$ in $W_j(r)$ is a singular torus (with
self-tangencies) which is a component of the level set $\{x\,|\,
\mu_{\Lambda_\infty}(x) = I_j\}$.
Since the pre-image of $T'_j$ in $\widetilde M_j$ has two distinct
components whose closures meet, there exists an element
of $\Gamma_j$ which does not commute with the image of
$\pi_1(T'_j)\to\Gamma_j$, but is contained in the image of
$\pi_1(\overline{T_j'})\to\Gamma_j$.  
It follows that 
\begin{claim}\label{schnauzer}
the image of $\pi_1(\overline{T_j'})$ in $\pi_1(W_j(r))$ is not carried
(up to conjugacy) by the frontier torus $\partial W_j(r)$.
\end{claim}

For the proof of (\ref{mr bluster}), we consider any subsequence
$(M_{j_n})$ of $(M_j)$ such that the manifolds $X_n =
f_{j_n}^{-1}(\overline{T'_{j_n}})$, which are contained in the compact
set $ W_\infty(r_1) $, converge in the Hausdorff topology.  The limit
is then a compact connected subset $X_\infty$ of $ W_\infty(r_1) $
which contains the frontier of $W_\infty(r_1)$.  Moreover it follows
from \ref{lisa} that $X_\infty = \{ x \,|\, \mu_{\Lambda_\infty}(x)
\le I \}$, where $I =\lim_{n\to\infty} I_{j_n}$.

Next we will show that $I = I_\infty$, and hence that
$X_\infty=\overline{T'_\infty}$.  If $I < I_\infty$ then $X_\infty$ is
a submanifold contained in the interior of $T'_\infty$ and hence, for
small $\epsilon>0$, the image of $\pi_1(N(\epsilon,X_\infty))$ is
carried by $\partial W_\infty(r)$.  However, for sufficiently large
$j$, we have that $X_j\subset N(\epsilon,X_\infty)$ and hence that the
image of $\pi_1(X_j)$ is carried by $\partial W_\infty(r)$.  Since
$f_j$ is a homeomorphism and $f_j(X_j) = T'_j$, this contradicts
\ref{schnauzer}.  If $I > I_\infty$ then, for sufficiently large $j$,
$f_j(\overline{T'_\infty})$ is contained in the interior of $T'_j$.
Thus the image of $\pi_1(f_j(\overline{T'_\infty}))$ is carried by
$\partial W_j(r)$, which again contradicts \ref{schnauzer}.

Since the manifolds $X_n$ are converging in the
Hausdorff topology to $\overline{T'_\infty}$  we have
$\vol X_n\to\vol T'_\infty$.  Moreover, since the Lipschitz
constants of the homeomorphisms $f_{j_n}$ are converging to $1$,
we also have $\vol T'_{j_n} \to \vol T'_\infty$
and $\vol(W_{j_n}(r_1) - T'_{j_n}) \to \vol(W_\infty(r_1) - T'_\infty)$ 
as $n\to\infty$.
Thus 
$$\vol (M_{j_n} - T_{j_n}) = \vol( W_{j_n}(r_1) - T'_{j_n}) \to 
\vol (W_\infty(r_1) - T'_\infty) = \vol( M_\infty - T_\infty) .$$
In view of conclusion (\ref{flubadub}) it follows that
$\lim_{n\to\infty} \vol T_{j_n} = \vol T_\infty$. 

Thus we have shown that for any subsequence $(M_{j_n})$ of $(M_j)$
such that $(X_n)$ converges in the Hausdorff topology, we have
$\lim_{n\to\infty} \vol T_{j_n} = \vol T_\infty$. Since every
subsequence of $(M_j)$ contains a subsequence
$(M_{j_n})$ of this type, conclusion (\ref{mr bluster}) follows.

Assume, for the proof of conclusion (\ref{dilly-dally}), that
$\lambda$ is a Margulis number for each $M_j$.  It suffices to show
that $\lambda - 2\delta$ is a Margulis number for $M_\infty$ for any
given $\delta$ with $ 0<\delta\le\delta_0$.

For a sufficiently large $r>r_1$, the frontier torus of $W_\infty(r)$
is contained in the $\lambda$-thin part of $M_\infty$. We
fix a number $r=r_2$ with this property, and denote the frontier torus
of $W_\infty(r_2)$ by $F$.  We consider a positive number
$\delta$ which is small enough to ensure that $F$ is contained in the
($\lambda-2\delta$ )-thin part of $M_\infty$, which we shall denote by
$\Theta$.

By \ref{lisa} and the compactness of $W_{\infty}$ we may fix $K\in\NN$
so that $\inj_{M_\infty}(x) \ge \inj_{M_K}(f_K(x)) - \delta$ for
all $x\in W_\infty(r_2)$. Hence if $\Psi$ denotes the $\lambda$-thin
part of $M_K$, we have $f_K(\Theta\cap W_\infty(r_2))\subset\Psi$, so
that $\Theta$ is contained in the set $Q=f_K^{-1}(\Psi\cap
W_K(r_2))\cup(M_\infty-W_\infty(r_2))\subset M_\infty$.

We are required to show that $\lambda - 2\delta$ is a Margulis number
for $M_\infty$.  According to the discussion in \ref{cylinder
  notions}, this is equivalent to showing that if $\theta$ is a
component of $\Theta$, then $\pi_1(\theta)$ has abelian image in
$\Gamma \pi_1(M_\infty)$. Hence it suffices to show that every
component of $Q$ has abelian fundamental group.

We have seen that $V_+=\overline{M_\infty-W_\infty(r_2)}$ is homeomorphic to
$S^1\times S^1\times[1,\infty)$, so that in particular $V_+$ is
connected. We have $F\subset\Theta\cap W_\infty(r_2)\subset
f_K^{-1}(\Psi\cap W_K(r_2))$. If $V_-$ denotes the component of
$f_K^{-1}(\Psi\cap W_K(r_2))$ containing $F$, then $V_+\cup V_-$ is one
component of $Q$, and every other component has the form
$f_K^{-1}(\psi)$, where $\psi$ is a component of $\Psi$ contained in
$\inter W_K(r_2)$. Any component of the latter type is homeomorphic to
a component of $\Psi$, and is therefore a solid torus since $\lambda$
is a Margulis number for $M_K$. In particular such a component has an
abelian fundamental group.

To describe the topology of $V_+\cup V_-$, we consider the component
$\psi_0$ of $\Psi$ containing $f_K(V_-)$. Since $\lambda$ is a
Margulis number for $M_K$, the set $\psi_0$ is a tube about some
closed geodesic $C_0$ in $M_K$. But $f_K(F)\subset\psi_0$ is the
frontier of $M_K-W_K(r_2)$, which we have seen is a topological
regular neighborhood of the core geodesic $C$ in $M_K$. We must
therefore have $C=C_0$. As $M_K-W_K(r_2)$ and $\psi_0$ are open
regular neighborhoods of $C_0$, with $M_K-W_K(r_2)\subset\psi_0$, the
set $\psi_0\cap W_K(r_2)=f_K(V_-)$ is homeomorphic to $S^1\times
S^1\times(0,1]$, and hence so is $V_-$. Since $V_+\cong S^1\times
S^1\times[1,\infty)$, and $\partial V_+=\partial V_-$, it follows that
$V_+\cup V_-\cong S^1\times S^1\times(0,\infty)$; in particular,
$\pi_1(V_+\cup V_-)$ is abelian.
\EndProof

The proof of the following result combines Propositions \ref{popeye}
and \ref{tube convergence}.

\Theorem\label{geometric fact} Let $M$ be a complete, finite-volume,
orientable hyperbolic manifold having exactly one cusp, and let
$k\ge3$ be an integer. Suppose that there is an infinite sequence of
closed manifolds $(M_j)_{j\in\NN}$, obtained by distinct Dehn fillings
of the compact core of $M$, such that $\pi_1(M_j)$ is $k$-free for
every $j\ge0$.  Let $H$ denote the maximal standard neighborhood
(see \ref{distance notions}) of the cusp of
$M$. Suppose that $\vol H<2\pi$. Then there exist a real number $\beta
$ with $1<\beta<2 $ and a point $Q\in \Mthick(\log3)$ such that
\Conclusions
\item $\vol H\ge \pi \beta $, and
\item $\dist(Q,\bar H)
  \ge-\frac12\log\bigg(\displaystyle\frac{\beta-1}{k-2}\bigg).$
\EndConclusions

\EndTheorem

\Proof According to Theorem \ref{Dehn filling for the lumpen
  proletariat}, all but finitely many of the $M_j$ are hyperbolic Dehn
fillings of $M$; hence we may assume, after passing to a subsequence,
that all the $M_j$ are hyperbolic Dehn fillings.

Let $c_j$ denote the core geodesic of the Dehn filling $M_j$ and set
 $l_j=\length c_j$.  According to Proposition \ref{tube convergence}
we have $l_j\to0$ as $j\to\infty$.

Since $\pi_1(M_j)$ is in particular $2$-free for each $M_j$, the
number $\log3$ is a Margulis number for each $M_j$ (see \ref{log 3
works}).

After passing to a smaller subsequence we may assume that $l_j<\log3$
for every $j$.  Let $T_j$ denote the maximal (embedded) tube about
$c_j $, let $R_j$ denote its radius, and set
$$R_j '=\frac12\arccosh\big(\cosh ( 2 R_j) \cosh\big(\frac{ l_j}2\big)\big).$$
Then according to Proposition \ref{popeye} there exist points $P_j \in
c_j $ and $Q_j\in M_j$ such that $\rho=\dist_M(P_j ,Q_j)$ satisfies
\Equation\label{olive oyl}\frac{k-2}{1+e^{2\rho}}+
\frac{1}{1+e^{l_j}}+\frac{1}{1+e^{2R_j' }} \le\frac12.\EndEquation
Furthermore, it follows from Proposition
\ref{popeye} that we may assume that our subsequence has been chosen
so that $Q_j$ lies in the $(\log3)$-thick part of $M_j$ for every $j$.

From the definition of the $R_j'$ we see that $R_j\le R_j'\le
R_j+( l_j/4)$. Since $l_j\to0$, it follows that $R_j-R_j'\to0$ as
$j\to\infty$. 

From (\ref{olive oyl}) it follows in particular that
$$\frac{1}{1+e^{l_j}}+\frac{1}{1+e^{2R_j' }}
\le\frac12.$$
Let us set 
$$r_j=\frac12\log\bigg(\frac{e^{l_j}+3}{e^{l_j}-1}\bigg),$$
so that
\Equation\label{reggie}\frac{1}{1+e^{l_j}}+\frac{1}{1+e^{2r_j }}
=\frac12.\EndEquation
Then we have $R_j'\ge r_j$ for every $j$. In particular,
$R_j'\to\infty$, and hence $R_j\to\infty$, as $j\to\infty$. 

We set $h_j'=R_j'-r_j\ge0$, and $h_j=R_j-r_j$. Note that
$h_j'-h_j\to0$ as $j\to\infty$.

According to Proposition \ref{tube convergence} we have
\Equation\label{wimpy}
\lim_{j\to\infty}\vol T_j = \vol H. 
\EndEquation 
On the other hand, we have $\vol T_j=\pi l_j\sinh^2 R_j$.  Since $R_j$
tends to infinity with $j$, it follows that
\Equation\label{swee'pea}
\frac{\vol T_j}{l_je^{2R_j}}\to\frac\pi4.
\EndEquation
We have $R_j=r_j+h_j$.  From the definition of $r_j$, and the fact
that $l_j\to0$, we find that $e^{2r_j}l_j\to 4$.  Combining these
observations with (\ref{wimpy}) and (\ref{swee'pea}), we deduce that
$\lim_{j\to\infty}\pi e^{2h_j}$ exists and that
$$\lim_{j\to\infty}\pi e^{2h_j}=\lim_{j\to\infty}\vol T_j  = \vol(H).$$
Hence $\alpha=\lim_{j\to\infty}h_j$ exists, and if we set
$\beta=e^{2\alpha}$ we have
\Equation\label{mehitabel}\vol H = \beta\pi.\EndEquation

Since $h_j'-h_j\to 0$ we have $h_j'\to\alpha$ as $j\to\infty$. In
particular, since $h_j'\ge0$, we have $\alpha\ge0$ and hence
\Equation\label{archy}\beta\ge1.\EndEquation
The hypothesis $\vol H<2\pi$, with (\ref{mehitabel}), gives
\Equation\label{boss}\beta<2.\EndEquation

By comparing (\ref{olive oyl}) and (\ref{reggie}) we find that 
$$\frac{k-2}{1+e^{2\rho_j}}+\frac1{1+e^{2R'_j}}\le\frac1{1+e^{2r_j}},$$
which gives
\Equation\label{veronica}
(k-2)\frac{1+e^{2R'_j}}{1+e^{2\rho_j}}
\le\frac{1+e^{2R'_j}}{1+e^{2r_j}}-1.
\EndEquation 
As $j\to\infty$, since $R'_j-r_j=h'_j\to\alpha$, the right hand side
of (\ref{veronica}) approaches $e^{2\alpha}-1=\beta-1$. Hence
$$(k-2)\limsup e^{2(R_j'-\rho_j)}\le\beta-1.$$
Since $R_j\le R_j'$ it follows that
$$(k-2)\limsup e^{2(R_j-\rho_j)}\le\beta-1.$$
Recalling that $\beta<2$ by (\ref{boss}), we find that
$$0< \liminf(\rho_j-R_j)\le\infty.$$  
Hence after passing to a subsequence we may assume that 
$\rho_j-R_j\to\gamma$, where $\gamma=+\infty$ if $\beta=1$ and
\Equation\label{dagwood}\gamma\ge-\frac12\log\bigg(\frac{\beta-1}{k-2}\bigg)
\qquad{\rm if}\qquad \beta>1.\EndEquation
In particular, since $\beta < 2$ and $k \ge 3$, we have $\gamma > 0$.

For every $j$ set $d_j=\dist(Q_j,\bar T_j)$ (so that, a priori, we
have $d_j\ge0$, with equality if and only if $Q_j\in\bar T_j)$. Since
$T_j=N(R_j,c_j)$, and $c_j$ has diameter at most $l_j/2$, we have
$\rho_j=\dist(P_j,Q_j)\le d_j+R_j+(l_j/2)$, and hence
\Equation\label{betty}d_j\ge\rho_j-R_j-(l_j/2).\EndEquation In
(\ref{betty}) we have $\rho_j-R_j\to\gamma$ and $l_j\to0$.  Setting
$\delta=\gamma/2$ if $0<\gamma<\infty$ and $\delta=1$ if
$\gamma=\infty$, we deduce that $d_j\ge\delta$ for large $j$,
i.e. $Q_j\in M_j-N(\delta,\bar T_j)$ for large $j$. In view of the
properties of the maps $f_j$ it then follows that $f_j(Q_j)\in M-H$
for every sufficiently large $j$. As $M-H$ is compact, we may assume
after passing to a subsequence that the sequence $(f_j(Q_j))$
converges to a point $Q\in M-H$. Furthermore, we have
$d_j\to\dist(Q,H)$ as $j\to\infty$. Taking limits in (\ref{betty}) we
deduce that \Equation\label{jughead}\gamma\le\dist(Q,H).\EndEquation
In particular $\gamma<\infty$, so that the definition of $\gamma$
gives $\beta>1$, and by (\ref{dagwood}) we have
$\gamma\ge-(1/2)\log((\beta-1)/(k-2))$. It now follows from
(\ref{jughead}) that
$$\dist(Q,\bar H)\ge-\frac12\log\bigg(\frac{\beta-1}{k-2}\bigg).$$

This is conclusion (2) of the theorem. Conclusion (1) is included in
(\ref{mehitabel}). This completes the proof of the Theorem.
\EndProof

\section{Volume and density}\label{packing section}

\Number\label{numbers}
We shall define functions $B(r)$, $a_2(r)$, $h_2(r)$, $a_3(r)$,
$h_3(r)$, $\beta(r)$, $\tau(r)$ and $d(r)$ for $r>0$. Each of these
will be defined by an analytic formula accompanied by a geometric
interpretation. For $n=2,3$ we let $\Delta_n(r)$ denote a regular
hyperbolic simplex in $\HH^n$ with sides of length $2r$.

\Bullets
\item $B(r)=\pi(\sinh(2r)-2r)$ is the volume of a ball of radius $r$
  in $\HH^3$;
\item $a_2(r)=\arccosh(\cosh(2r)/\cosh(r)$ is the altitude of
  $\Delta_2(r)$;
\item $h_2(r)=\arctanh((\cosh a_2(r)\cosh(r)-1)/\sinh a_2(r)\cosh r)$
  is the distance from a vertex of $\Delta_2(r)$ to its barycenter;
\item $a_3(r)=\arccosh(\cosh(2r)/\cosh(h_2(r)))$ is the altitude of
  $\Delta_3(r)$;
\item $h_3(r)=\arctanh((\cosh a_3(r)\cosh h_2(r)-1)/ (\sinh
  a_3(r)\cosh h_2(r)))$ is the distance from a vertex of $\Delta_3(r)$
  to its barycenter;
\item $\beta(r)=\arcsec(\sech(2r)+2)$ is a dihedral angle of
  $\Delta_3(r)$;
\item $\tau(r)=3\int_{\beta(r)}^{\arcsec3}\arcsech((\sec t)-2)dt$ is
  the volume of $\Delta_3(r)$;
\item $d(r)=(3\beta(r)-\pi)(\sinh(2r)-2r)/\tau(r)$ is the ratio of
  $\vol(W\cap\Delta_3(r))$ to $\vol\Delta_3(r)$, where $W$ is the
  union of four mutually tangent balls of equal radii centered at the
  vertices of $\Delta_3(r)$.
\EndBullets

It is shown in \cite{bf} that $d(r)$ is a monotonically increasing
function for $r>0$. We set
$d(\infty)=\lim_{r\to\infty}d(r)=.853\ldots$
\EndNumber

If $P$ is a point in a complete hyperbolic manifold $M$, and if the
injectivity radius at $P$ is at least $r$ then it is immediate that
$\vol N(r,P) = B(r)$.  Less trivially, we have:

\Proposition\label{packing}Suppose $r$ is a positive number, that $P$
is a point in a hyperbolic manifold $M$, and that the injectivity
radius of $M$ at $P$ is at least $r$. Then  $\vol(N(h_3(r),P))\ge
B(r)/d(r)$.  \EndProposition

\Proof Set $J=N(r,P)$ and $W=N(h_3(r),P)$. Then $J$ is isometric to a
ball of radius $r$ in $\HH^3$, since the injectivity radius of $M$ at
$P$ is at least $r$. If we write $M=\HH^3/\Gamma$ where
$\Gamma\le\Isom_+(\HH^3)$ is discrete and torsion-free, and if
$q:\HH^3\to M$ denotes the quotient map, then each component of
$q^{-1}(J)$ is a hyperbolic ball of radius $r$. These components form
a sphere-packing in the sense of \cite{bor}. Let us fix any component
$\widetilde J$ of $q^{-1}(J)$, and let $\widetilde W\subset\HH^3$
denote the ball of radius $h_3(r)$ with the same center as $\widetilde
J$. The {Vorono\u\i} polytope $Z$ for $\widetilde J$, whose interior
consists of all the points which are nearer to the center of
$\widetilde J$ than to the center of any other sphere of the packing,
is a Dirichlet domain for $\Gamma$.  In the proof of the main result
of \cite{bor} (see \cite[p. 259]{bor}) it is shown that
$(\vol\widetilde J)/\vol(Z\cap\widetilde W)\leq d(r)$.  But
$\vol\widetilde J=B(r)$; and since $q$ maps $\widetilde W$ into $W$,
and is one-to-one on the interior of $Z$, we have
$\vol(Z\cap\widetilde W)\leq \vol W$. Hence $B(r)/(\vol W)\leq d(r)$.
\EndProof

\Proposition\label{horoball packing}Suppose that $H$ is a standard
neighborhood of a cusp in a finite-volume hyperbolic $3$-manifold $M$, and
set 
$$N=N(\log\sqrt{\scriptstyle\frac{3}{2}}\;,\bar H)\subset M.$$ 
Then
$$\frac{\vol H}{\vol N}\le d(\infty).$$
\EndProposition

\Proof This is essentially a restatement of \cite[Theorem 4]{bor}.  In
\cite[Section 6]{bor} it is remarked that the proof of the main
theorem (about sphere packings) applies without change to horoball
packings.  A stronger version of \cite[Theorem 4]{bor}, analogous to
the statement which was used in the previous proposition, follows from
this proof in the same way.  However, it must be formulated slightly
differently since $h_3(r)\to\infty$ as $r\to\infty$.  For this
reformulation we consider the function $k_3(r) \dot= h_3(r) - r$.
Geometrically, the value of $k_3(r)$ can be described as the distance
from the barycenter of $\Delta_3(r)$ to any of the four spheres of
radius $r$ centered at a vertex of $\Delta_3(r)$.  Note that if $P$ is
a vertex of $\Delta_3(r)$ and $B$ is a ball of radius $r$ centered at
$P$ then $N(h_3(r),P)=N(k_3(r),B)$.

If $\Delta$ is a regular ideal hyperbolic simplex then there is a
unique family of four ``standard'' horoballs which meet the sphere
$S^\infty$ at the vertices of $\Delta$ and are each tangent to the
other three.  By the barycenter of $\Delta$ we shall mean the unique
common fixed point of the 12 symmetries of $\Delta$, which permute the
standard horoballs.  It is clear that the quantity $k_3(\infty) =
\lim_{r\to\infty}k_3(r)$ is the distance from the barycenter of
$\Delta$ to any of these four horoballs.  Working in the upper
half-space model, with coordinates $(z,t)$, where $z\in \CC$ and $t\in
\RR_+$, it is straightforward to check, for the regular ideal hyperbolic
tetrahedron with vertices $(0,0)$, $(1,0)$, $(\frac{1+i\sqrt3}{2},0)$ and
$\infty$, that the barycenter is the point $(\frac{3+i\sqrt3}{6},
\sqrt\frac{2}{3})$ and that the standard horoball containing $\infty$
is bounded by the plane $t=1$.  Thus we have $k_3(\infty) =
\log\sqrt\frac{3}{2}$.

The components of the pre-image of $H$ in $\HH^3$ form a horoball
packing to which \cite[Theorem 4]{bor} applies, and a Ford domain $Z$
for one of the horoballs is a {Vorono\u\i} polyhedron for the packing.  The
theorem states that ${\vol H}/{\vol Z} \le d(\infty)$.  However, the
proof which was indicated for the theorem actually shows that ${\vol
  H}/{\vol(Z\cap N(k_3(\infty),H)} \le d(\infty)$.  Since $k_3(\infty)
=\log\sqrt\frac{3}{2}$, this is equivalent to the statement that
${\vol H}/{\vol N} \le d(\infty)$.
\EndProof

\section{The main estimate}\label{main estimate section}

This section contains the proof of Theorem \ref{main estimate},
which was stated in the Introduction.  

\Lemma\label{margulis}
Let $M$ be a complete, finite-volume, orientable hyperbolic manifold
having exactly one cusp. Suppose that there is an infinite sequence of
closed manifolds $(M_j)_{j\in\NN}$, obtained by distinct Dehn fillings
of the compact core of $M$, such that $\pi_1(M_j)$ is $2$-free for
every $j\ge0$. Then $\log3$ is a Margulis number for $M$.
\EndLemma

\Proof According to Theorem \ref{Dehn filling for the lumpen
  proletariat}, we may assume, after removing finitely many terms from
the sequence $(M_j)$, that the $M_j$ are all hyperbolic.  According to
\ref{log 3 works}, the number $\log3$ is a Margulis number for each
$M_j$, and it therefore follows from Proposition \ref{tube
  convergence} that $\log3$ is a Margulis number for $M$.
\EndProof

\Proof[Proof of Theorem \ref{main estimate}]
Let $H$
denote the maximal standard neighborhood  of the cusp of $M$. If $\vol
H\ge2\pi$ the conclusion is obvious, and so we shall assume  $\vol
H<2\pi$. According to Theorem \ref{geometric fact} we may fix a real number $\beta $ with $1<\beta<2 $ and a point $Q\in \Mthick(\log3)$ such
that
\Conditions
\item $\vol H\ge \pi \beta $, and
\item $\dist(Q,\bar H)\ge-(\log(\beta -1))/2$.
\EndConditions
Let us set $$N=N(\log\sqrt{\scriptstyle\frac{3}{2}},\bar H)\subset
M.$$ According to Proposition \ref{horoball packing}, we have $\vol
H/\vol N\le d(\infty)$, and by (1) it follows that
$$\vol N\ge \pi \beta/d(\infty) .$$
It follows from (2) that 
$$\dist(Q,\bar N)\ge-\frac{\log(\beta -1)}2-\frac
{\log(3/2)}2,$$
i.e.
$$\dist(Q,\bar N)\ge-\frac12\log\bigg(\frac32(\beta -1)\bigg).$$
Setting 
$$r_\beta=-\frac12\log\bigg(\frac32(\beta -1)\bigg)$$
and
$$N_1=N(r_\beta,Q)$$
we deduce that $N\cap N_1=\emptyset$, and hence that
\Equation\label{daffy}
\vol M\ge\frac{\pi \beta}{d(\infty)}+\vol N_1.
\EndEquation

We shall now distinguish several cases depending on the size of
$r_\beta$. To simplify the notation we shall set
$$C_1=\exp\bigg(h_3\bigg(\frac{\log3}2\bigg)\bigg)=1.931\ldots$$
and
$$C_2=\frac{B((\log3)/2)}{d((\log3)/2)}=0.929\ldots$$
\medskip

{\bf Case I:  $r_\beta \le 0$. } 

In this case $N_1=\emptyset$.  However, we have $\beta\ge5/3$, hence
$\pi \beta/d(\infty)>6$, so the conclusion follows from (\ref{daffy}).

\medskip

{\bf Case II:  $0<r_\beta\le(\log3)/2$.} 

In this case we have $5/3>\beta\ge11/9$.  Since $Q\in \Mthick(\log3)$,
the set $N_1$ is a hyperbolic ball of radius $r_\beta$. In the
notation of \ref{numbers} we have $\vol
N_1=B(r_\beta)=\pi(\sinh(2r_\beta)-2r_\beta)$, and (\ref{daffy})
becomes
$$\vol M\ge\pi\bigg(\frac{\beta}{d(\infty)}-\frac34(\beta-1)-\frac13(\beta-1)^{-1}
+\log\bigg(\frac32(\beta-1)\bigg)\bigg);$$
that is, we have $\vol M\ge f(\beta-1)$, where the function $f(x)$ is
defined for $ x>0$ by 
$$f(x)= \pi\bigg(\frac{ x+1}{d(\infty)}-\frac34x+\frac13x^{-1}
+\log\bigg(\frac32x\bigg)\bigg).$$
We have $f(x)\to\infty$ as $x$ approaches $0$ or
$\infty$, and $f'(x)=\pi q(1/x)$, where $q$ is the quadratic
polynomial
$$q(y)=-\frac13y^2+y+( \frac1{d_\infty}-\frac34).$$
As the only zero of $q(y)$ with $y>0$ is at
$y=(3/2)(1+\sqrt{4/(3d_\infty)})$, the least value of $f(x)$ on $(0,\infty)$ is
$$f\bigg(\frac2{3(1+\sqrt{4/(3d_\infty)})}\bigg)>\whatvolume.$$
In particular it follows that $\vol M>\whatvolume$ in this case.

\medskip

{\bf Case III:  $(\log3)/2 < r_\beta\le h_3((\log3)/2)$.}

In this case $11/9>\beta\ge1+\frac23C_1^{-2}=1.178\ldots$.  We have
$N_1\supset N((\log3)/2,Q)$ and, since $Q\in \Mthick(\log3)$, the set
$N((\log3)/2,Q)$ is a hyperbolic ball of radius $(\log3)/2$. Hence
$$\vol N_1\ge  B((\log3)/2) =0.737\ldots$$
We also have $\beta\ge1.178$, and so (\ref{daffy}) gives
$$\vol M\ge\frac{1.178\pi}{d(\infty)}+.737\ldots>5.07.$$

\medskip

{\bf Case IV:  $r_\beta > h_3((\log 3)/2)$.}

In this case we have $\beta<1+\frac23C_1^{-2}$.  We
set $$s_\beta=\frac12(r_\beta-h_3((\log3)/2).$$

Lemma \ref{margulis} guarantees that $\log3$ is a Margulis number for
$M$. Hence some component $X$ of the $(\log3)$-thin part of
$\Mthin(\log3)$ is a standard neighborhood of the cusp of $M$. Since
$H$ is the unique maximal standard neighborhood of the cusp, we have
$X\subset H$. It follows that the frontier $Y$ of $H$ is contained in
$\Mthick(\log3)$.

Since $\log3$ is a Margulis number for $M$, the set $\Mthick(\log3)$
is connected (see \ref{cylinder notions}). Define a continuous
function $\Delta :\Mthick(\log3)\to\RR$ by $\Delta
(x)=\dist_M(x,Q)$. Then $\Delta (Q)=0$. For any point $y\in
Y\subset\Mthick(\log3)$ we have
$\Delta(y)\ge\dist(Q,H)\ge-\log(\beta-1)/2> r_\beta$. By the
connectedness of $\Mthick(\log3)$ we have $[0,r_\beta ]\subset
\Delta(\Mthick(\log3))$.  Since the definition of $s_\beta$ and the
hypothesis of Case IV give that $0<2s_\beta+h_3((\log3)/2)<r_\beta$, we
may choose a point $R\in \Mthick(\log3)$ such that
$\Delta(R)=2s_\beta+h_3((\log3)/2)$, i.e.
$\dist(R,Q)=2s_\beta+h_3((\log3)/2)$. This implies on the one hand that
$N_2=N(2s_\beta,R)$ is disjoint from $W=N(h_3((\log3)/2),Q)$, and on
the other hand that
$$N_2\subset N(2s_\beta+h_3((\log3)/2),Q)=N(r_\beta,Q)=N_1.$$
Hence 
\Equation\label{porky}\vol N_1\ge\vol N_2+\vol W.\EndEquation

According to Proposition \ref{packing}, we have
\Equation\label{priscilla}\vol W\ge C_2.\EndEquation Combining
(\ref{daffy}), (\ref{porky}) and (\ref{priscilla}), we obtain
\Equation\label{donald}\vol M>C_2+\frac{\pi \beta}{d(\infty)}+\vol
N_2\EndEquation

We have two sub-cases.

\medskip

{\bf Case IV(A):  $ h_3((\log3)/2) < r_\beta\le h_3((\log3)/2)+\log3$. }

In this sub-case we have $1+\frac23C_1^{-2} >\beta\ge
1+\frac29C_1^{-2} =1.05\ldots$.  Also, since $r_\beta\le
h_3((\log3)/2)+\log3$, we have $s_\beta\le(\log3)/2$.

Since $R\in \Mthick(\log3)$, the set $N_2=N(s_\beta,R)$ contains a
hyperbolic ball of radius $s_\beta$. Hence $\vol N_2\ge B(s_\beta)$,
and (\ref{donald}) gives

\Equation\label{daisy}
\vol M>C_2+\frac{\pi \beta}{d(\infty)}+B(s_\beta).
\EndEquation

Recalling that $B(s_\beta)=\pi(\sinh(2s_\beta)-2s_\beta)$, that
$s_\beta= (1/2)(r_\beta-h_3((\log3)/2)$, and that $
r_\beta=-(1/2)\log\big(\frac32(\beta -1)\big)$, we obtain the
expression
$$B(s_\beta)=\frac\pi2\bigg(C_1^{-1}\big(\frac32(\beta-1)\big)^{-1/2}-C_1\big(\frac32(\beta-1)\big)^{1/2}+\log\big(\frac32(\beta-1)\big)+2\log
C_1\bigg).$$
We may therefore rewrite
(\ref{daisy}) as $\vol M\ge g(\beta-1)$, where $g(x)$ is defined for
$x>0$ by
$$g(x)=\frac\pi2\bigg(C_1^{-1}\big(\frac32x\big)^{-1/2}-C_1\big(\frac32x\big)^{1/2}+\log\big(\frac32x\big)+2\log C_1+\frac2{d(\infty)}(x+1)\bigg)+C_2.$$

Note that for $x>0$ we have 
$$g'(x)=\frac\pi2c\bigg(\sqrt{\frac3{2x}}\bigg),$$
where $c$ denotes the cubic polynomial 
$$c(y)=-\frac34C_1^{-1}y^3+\frac23y^2-\frac34C_1y+\frac2{d_\infty}.$$
Writing
$$c(y)=-\frac34C_1^{-1}y^2(y-A)-\frac34C_1(y-B),$$
where $A=(8/9)C_1=1.7\ldots$ and $B=8/(3d_\infty C_1)=1.6\ldots$,
makes it evident that $c(y)<0$ for $y>2$. Hence $g(x)$ decreases
monotonically on the interval $0<x<3/8$. Since we are in Case IV we
have $\beta-1<(2/3)C_1^{-2}=0.178\ldots<3/8$ and hence
$$\vol M\ge g(\beta-1)>g\big(\frac23C_1^{-2}\big)=
\frac\pi{d(\infty)}(1+\frac23C_1^{-2})+C_2=5.26\ldots
$$

\medskip

{\bf Case IV(B):  $r_\beta>h_3((\log3)/2)+\log3$}

In this sub-case we have $\beta<1+\frac29C_1^{-2}$, and
$s_\beta>(\log3)/2$. Since $R\in \Mthick(\log3)$, the set
$N_2=N(s_\beta,R)$ contains a hyperbolic ball of radius
$(\log3)/2$. Hence $\vol N_2\ge B((\log3)/2)$, and (\ref{donald})
gives
$$\vol M>C_2+\frac{\pi
  \beta}{d(\infty)}+B((\log3)/2)>C_2+\frac{\pi
  }{d(\infty)}+B((\log3)/2)=5.34\ldots$$
Thus the desired lower bound is established in all cases.
\EndProof

\section{Relaxing the restriction on surface subgroups}\label{singular applications}

In this section we prove Proposition \ref{lots of H_1}, which was
discussed in the introduction to the paper.

\Proposition\label{lots of H_1}
Let $M$ be a complete, finite-volume, orientable hyperbolic manifold
having exactly one cusp, such that 
 $\dim_{\ZZ_2}H_1(M;\ZZ_2)\ge 7$. 
Then either
\Alternatives
\item $\vol M>\whatvolume,$ or
\item $M$ contains a genus-$2$ connected incompressible surface.
\EndAlternatives 
\EndProposition

\Proof
Let us fix a basis $\{\lambda,\mu\}$ of $H_1(\partial N;\ZZ)$ such
that $\lambda$ lies in the kernel of the natural homomorphism
$H_1(\partial N;\ZZ)\to H_1(N;\ZZ_2)$. For every integer $n$, let
$\alpha_n$ denote the slope defined by the primitive class
$\lambda+2n\mu\in H_1(\partial N;\ZZ)$, and let $M_n$ denote the
closed manifold obtained from $N$ by the Dehn filling corresponding to
$\alpha_n$.  If $i:N\to M_n$ denotes the inclusion homomorphism, then
$i_*:H_1(N;\ZZ_2)\to H_1(M_n;\ZZ_2)$ is an isomorphism for every
$n$. In particular we have $\dim_{\ZZ_2}H_1(M_n;\ZZ_2)\ge7$.

We first consider the case in which there is an infinite sequence of
distinct integers $(n_j)_{j\in\NN}$ such that $\pi_1(M_{n_j})$ is
$3$-free for every $j\ge0$. In this case it follows immediately from
Theorem \ref{main estimate} that $\vol M>\whatvolume$. This is conclusion
(a) of the present proposition.

We shall therefore assume for the rest of the proof that there is an
integer $n_0>0$ such that $\pi_1(M_{n})$ fails to be $3$-free whenever
$|n|\ge n_0$. Since $\dim_{\ZZ_2}H_1(M_n;\ZZ_2)$ is in particular
$\ge5$, it follows from Proposition \ref{variant} that $\pi_1(M_n)$
has a subgroup isomorphic to a genus-$2$ surface group whenever
$|n|\ge n_0$. Since $\dim_{\ZZ_2}H_1(M_n;\ZZ_2)\ge7$, it then follows
from the case $g=2$ of \cite[Theorem \topeleven]{last} that
for each $n$ with $|n|\ge n_0$, the manifold $M_n$ contains a
connected incompressible closed surface of genus $2$.

The main theorem of \cite{hatcher} implies that there are only
finitely many boundary slopes (see \ref{hatcher}) in $\partial N$.
Hence there is an integer $n_1\ge n_0$ such that $\alpha_n$ is not a
boundary slope. We may write $M_{n_1}=N\cup T$, where $T$ is a solid
torus with $T\cap N =\partial T= \partial N $, and the meridian disks
of $T$ represent the slope $\alpha_{n_1}$.

Let $Y\subset M_{n_1}$ be a surface such that
\begin{enumerate}[(1)]
\item\label{scrooge} $Y$ is orientable, incompressible, connected
and closed, and $\genus(Y)=2$.
\end{enumerate}

We may choose $Y$ so that

\begin{enumerate}[(1)]
\setcounter{enumi}{1}
\item\label{goofy} each component of $Y\cap T$ is a meridian disk in
  $T$.
\end{enumerate}

We let $m$ denote the number of components of $Y\cap T$, and we
suppose $Y$ to be chosen, among all surfaces in $M$ for which
(\ref{scrooge}) and (\ref{goofy}) hold, so that $m$ is as small as
possible.

Suppose that $m>0$. Set $F=Y\cap N$, so that $F$ is a compact
orientable surface of genus $2$, properly embedded in $N $,
with $\partial F\ne\emptyset$.  Since $\partial N$ is a torus, $F$
cannot be boundary-parallel.  If the inclusion homomorphism
$\pi_1(F)\to\pi_1(N)$ were injective then $\alpha_{n_1}$ would be a
boundary slope, in contradiction to our choice of $n_1$.  On the
other hand, if $\pi_1(F)\to\pi_1(N )$ has a non-trivial kernel, it
follows from Dehn's lemma and the loop theorem that there is a disk
$D_1\subset N $ such that $D_1\cap F =\partial D_1$, and such that
$\partial D_1$ does not bound a disk in $F$. However, since $Y$ is
incompressible in $M_{n_1}$, the curve $\partial D_1$ does bound a
disk $D\subset Y$, which must contain at least one component of $Y\cap
T$. If we set $Y_1=(Y-D)\cup D_1$, then (\ref{scrooge}) and
(\ref{goofy}) still hold when $Y$ is replaced by $Y_1$, but $Y_1\cap
T$ has at most $m-1$ components. This contradicts the minimality of
$m$. Hence we must have $m=0$. Thus $Y\subset N$ is 
a genus-$2$ connected incompressible surface, and conclusion (b) of
the proposition holds.
\EndProof

As we mentioned in the introduction, in \cite{cds} it will be shown
that the hypothesis of Proposition \ref{lots of H_1} implies
conclusion (a). The proof of this stronger result uses Proposition
\ref{lots of H_1}.

\bibliographystyle{plain}
%\bibliography{howfree_v3.bib}

\end{document}